\documentclass[]{amsart}
\pdfoutput=1

\usepackage{latexsym}   
\usepackage{amssymb}    
\usepackage[usenames,dvipsnames]{xcolor} 
\usepackage[colorlinks, 
            linkcolor=MidnightBlue, 
            citecolor=MidnightBlue, 
            urlcolor =BlueViolet,
            breaklinks=true
           ]{hyperref}
\usepackage{booktabs}
\usepackage[all]{xy}
\usepackage{multirow}
\usepackage{colonequals}
\usepackage{xspace}

\newtheorem{theorem}{Theorem}[section]
\newtheorem{lemma}[theorem]{Lemma}
\newtheorem{proposition}[theorem]{Proposition}
\newtheorem{corollary}[theorem]{Corollary}
\newtheorem{conjecture}[theorem]{Conjecture}
\newtheorem{heuristic}[theorem]{Heuristic}

\newtheorem{corollarynote}[theorem]
{Corollary\footnote{Technically, this result is not a corollary but rather a \emph{porism} --- that is,
                    something that follows easily from the \emph{proof} of an earlier result, but not from its statement.}}

\theoremstyle{definition}
\newtheorem{remark}[theorem]{Remark}
\newtheorem{notation}[theorem]{Notation}
\newtheorem{definition}[theorem]{Definition}
\newtheorem{example}[theorem]{Example}

\numberwithin{equation}{section}

\DeclareMathOperator{\Aut}{Aut}

\DeclareMathOperator{\End}{End}

\DeclareMathOperator{\Hom}{Hom}
\DeclareMathOperator{\Jac}{Jac}

\newcommand \bA {{\mathbb A}}

\newcommand \FF {{\mathbb F}}
\newcommand \Fq {\FF_q}
\newcommand \PP {{\mathbb P}^1}
\newcommand \QQ {{\mathbb Q}}
\newcommand \ZZ {{\mathbb Z}}

\newcommand \CA {{\mathcal A}}
\newcommand \CF {{\mathcal F}}

\newcommand{\babbage}{\texttt{babbage}}
\newcommand{\hensel}{\texttt{hensel}}
\newcommand{\winifred}{\texttt{winifred}}

\newcommand{\abar}{\bar{a}}
\newcommand{\bbar}{\bar{b}}
\newcommand{\cbar}{\bar{c}}
\newcommand{\dbar}{\bar{d}}

\newcommand{\Dbar}{\bar{D}}
\newcommand{\Fbar}{\overline{\FF}}
\newcommand{\Kbar}{\bar{K}}
\newcommand{\Qbar}{\overline{\QQ}}

\newcommand{\Ctilde}{\tilde{C}}

\newcommand{\phihat}{\hat{\varphi}}

\newcommand{\psihat}{\hat{\psi}}

\newcommand{\bad}{\textup{bad}}
\newcommand{\sing}{\textup{sing}}
\newcommand{\toobig}{\textup{big}}
\newcommand{\pz}{\phantom{0}}   
\newcommand{\eps}{\varepsilon}

\newcommand{\onea}{\textup{1a}}
\newcommand{\oneb}{\textup{1b}}
\newcommand{\twoa}{\textup{2a}}
\newcommand{\twob}{\textup{2b}}
\newcommand{\threea}{\textup{3a}}
\newcommand{\threeb}{\textup{3b}}
\newcommand{\threec}{\textup{3c}}

\newcommand{\dsixcurve}{genus-$2$ $D_6$ curve\xspace}
\newcommand{\dsixcurves}{genus-$2$ $D_6$ curves\xspace}
\newcommand{\dsixfamily}{genus-$2$ $D_6$ family\xspace}

\newcommand{\mybar}[1]{
  \mathchoice
  {#1\llap{$\overline{\phantom{\displaystyle\rm#1}}$}}
  {#1\llap{$\overline{\phantom{\textstyle\rm#1}}$}}
  {#1\llap{$\overline{\phantom{\scriptstyle\rm#1}}$}}
  {#1\llap{$\overline{\phantom{\scriptscriptstyle\rm#1}}$}}
}  
\renewcommand{\bar}{\mybar}
\renewcommand{\tilde}{\widetilde}
\renewcommand{\hat}{\widehat}

\makeatletter
\@namedef{subjclassname@2020}{%
  \textup{2020} Mathematics Subject Classification}
\makeatother

\makeatletter

\def\@marginparreset{\marginparstyle}
\def\marginparstyle{\SMALL\normalfont\raggedright\openup-2pt }
\long\def \@savemarbox #1#2{%
\global\setbox #1%
  \vtop{%
    \hsize\marginparwidth
    \@parboxrestore
    \reset@font
    \@setnobreak
    \@setminipage
    \@marginparreset
    #2%
    \par
    \global\@minipagefalse
    }%
}

\DeclareRobustCommand\marginparhere[2][0pt]{%
\ifhmode\unskip\fi
\ifmmode\ssty\mathclose{\fi
\rlap{\hskip\marginparsep\smash{%
  \vtop to 0pt{\marginparstyle \hsize\marginparwidth 
    \leftskip=#1 \rightskip=-#1 plus20pt \noindent#2\vss}}}%
\ifmmode}\fi
}

\DeclareRobustCommand{\redden}{\@ifnextchar*
{\@latex@error{{redden*} is only an environment}}
{\@ifnextchar[{\r@dden}{\r@dd@n}}}
\def\r@dden[#1]{\ifmmode\@mperr{math}\else\ifinner\@mperr{inner}%
\else\leavevmode\marginpar{\leavevmode#1\endgraf}\fi\fi\r@dd@n}
\def\r@dd@n{\def\reserved@a{redden}
\ifx\@currenvir\reserved@a\redd@n\bgroup\ignorespaces
\else\expandafter\redd@n\fi}
\def\endredden{\unskip\egroup}
\def\@mperr#1{\@latex@warning{redden in #1 mode: no marginpar possible}}
\def\redd@n#1{\leavevmode{\color{red}#1}}

\makeatother

\begin{document}

\title[Doubly isogenous curves of genus two]
      {Doubly isogenous curves of genus two\\ with a rational action of \texorpdfstring{$D_6$}{D6}}

\date{\today}

\author[Booher]{Jeremy Booher}
\address[Booher]{Department of Mathematics,
1400 Stadium Rd.,
University of Florida,
Gainesville, FL 32611, USA}
\email{jeremybooher@ufl.edu}
\urladdr{https://people.clas.ufl.edu/jeremybooher}

\author[Howe]{Everett W. Howe}
\address[Howe]{Unaffiliated mathematician, 
         San Diego, CA 92104, USA}
\email{however@alumni.caltech.edu}
\urladdr{https://ewhowe.com}

\author[Sutherland]{Andrew V. Sutherland}
\address[Sutherland]{Department of Mathematics, 
         Massachusetts Institute of Technology, 
         77 Massachusetts Avenue, 
         Cambridge, Massachusetts 02139, USA}
\email{drew@math.mit.edu}
\urladdr{https://math.mit.edu/~drew}

\author[Voloch]{Jos\'e Felipe Voloch}
\address[Voloch]{School of Mathematics and Statistics, 
         University of Canterbury, Private Bag 4800,
         Christchurch 8140, New Zealand}
\email{felipe.voloch@canterbury.ac.nz}
\urladdr{https://www.math.canterbury.ac.nz/~f.voloch}

\thanks{Booher and Voloch were partly supported by the Marsden Fund administered by the Royal Society of New Zealand, Howe was partly supported by the Erskine Fellowship from the University of Canterbury and Sutherland was supported by Simons Foundation grant 550033.}



\begin{abstract}
Let $C$ and $C'$ be curves over a finite field $K$, provided with embeddings
$\eps$ and $\eps'$ into their Jacobian varieties. Let $D\to C$ and $D'\to C'$
be the pullbacks (via these embeddings) of the multiplication-by-$2$ maps on the
Jacobians. We say that $(C,\eps)$ and $(C',\eps')$ are \emph{doubly isogenous}
if $\Jac C$ and $\Jac C'$ are isogenous over $K$ and $\Jac D$ and $\Jac D'$
are isogenous over~$K$.  When we restrict attention to the case where $C$ and 
$C'$ are curves of genus~$2$ whose groups of $K$-rational automorphisms are 
isomorphic to the dihedral group $D_6$ of order~$12$, we find many more doubly 
isogenous pairs than one would expect from reasonable heuristics. 

Our analysis of this overabundance of doubly isogenous curves over finite fields
leads to the construction of a pair of doubly isogenous curves over a number 
field. That such a global example exists seems extremely surprising. We show
that the Zilber--Pink conjecture implies that there can only be finitely many
such examples. When we exclude reductions of this pair of global curves in our
counts, we find that the data for the remaining curves is consistent with our
original heuristic.

Computationally, we find that if $C$ and $C'$ are doubly isogenous curves in
our family of \dsixcurves, then $C$ and $C'$ have naturally defined unramified
abelian Galois covers whose Pryms are \emph{not} isogenous to one another; in
practice, it is enough to consider certain covers with Galois groups of exponent $3$ 
and~$4$.

We discuss how our family of curves can potentially be used to obtain a
deterministic polynomial-time algorithm to factor univariate polynomials over
finite fields via an argument of Kayal and Poonen.
\end{abstract}

\maketitle


\section{Introduction}
\label{sec:introduction}

Let $C$ be a (smooth, geometrically irreducible, projective) algebraic curve of genus $g>1$ over a finite field $k$ of odd characteristic. Choose an embedding 
$\eps\colon C \to J$ from $C$ to its Jacobian $J$ by fixing a divisor of degree one. The pullback of the isogeny $[2]\colon J \to J$ by $\eps$ is a curve $C^{[2]}$ that is an \'{e}tale cover of~$C$.
(This depends on the embedding $\eps$; as this will be clear from context we do not include it in the notation.)
Sutherland and Voloch raised the question (expressed in different terms) of whether the $L$-function of $C^{[2]}$ characterizes $C$ if $g>2$ \cite[Question~2.3]{SutherlandVoloch2019}. They also provided a heuristic argument and numerical evidence for a positive answer. The same heuristics pointed to a negative answer when $g=2$ and, indeed, an example was provided there. More formally: 

\begin{definition}
\label{def:doubly isogenous}
Two (smooth, geometrically irreducible, and projective) algebraic curves $(C,\eps)$ and $(C',\eps')$ of genus $g>1$ over a field $k$ of odd characteristic, together with embeddings into their respective Jacobians,
are \emph{doubly isogenous} if their Jacobians are isogenous and, additionally, the Jacobians of the curves $C^{[2]}$ and $C^{'[2]}$ defined above are also isogenous. 
\end{definition}

In the present paper we study a one-parameter family of curves of genus~$2$: the ``\dsixcurves,'' whose automorphism groups contain the dihedral group $D_6$ of order~$12$. In Section~\ref{sec:family} we review basic results about this family,
including the fact that the Jacobians of such curves are isogenous to the square of an elliptic curve.
While searching for examples of doubly isogenous genus-$2$ curves, we found many more examples involving curves from this family than we na\"{\i}vely expected.
Our goal is to explain this observation and to use information from additional covers to distinguish curves in this family from one another.

For \dsixcurves $C$, we study the decomposition of the Jacobians of $C^{[2]}$ as a product of simple abelian varieties in Section~\ref{sec:covers2}, and in Section~\ref{sec:2coverheuristics} we formulate a refined heuristic for the number of doubly isogenous curves in this family based on the symmetries of this decomposition.
However, numerical data over finite fields presented in Section~\ref{sec:2coverheuristics} does not match our refined heuristic. 
In Section \ref{sec:extraordinary} we explain this discrepancy by showing that there exists an example of a pair of doubly isogenous curves in characteristic zero that influences the statistics we observe over finite fields.

\begin{theorem}[cf. Theorem~\ref{T:extraordinary}]
Let $K = \QQ(\sqrt{29})$.  Let $L$ be the $S_3$ extension of $K$ obtained by adjoining the roots of
$x^3 + x^2 + 2$ and $x^2 + 1$.  There exist \dsixcurves $C_1$ and $C_2$ defined over~$K$, with Weierstrass points defined over~$L$, such that $C_1$ and $C_2$ are doubly isogenous over $L$ \textup(when embedded into their Jacobians using a Weierstrass point\textup).  
\end{theorem}

We refer to $C_1$ and $C_2$ as \emph{the extraordinary curves}, because their existence is extremely surprising --- see Remark~\ref{remark:simultaneous}.

Over some finite fields, certain pairs of \dsixcurves that are twists of one another are also automatically doubly isogenous (see Theorem~\ref{T:easytwists}).  After excluding the reductions of $C_1$ and $C_2$ and the doubly isogenous twists from Theorem~\ref{T:easytwists} from our data set, we find that
the number of remaining doubly isogenous curves in the \dsixfamily agrees with our refined heuristic, as can be seen in Table~\ref{table:noextraordinary}.  See Sections \ref{sec:2coverheuristics} and \ref{sec:extraordinary} for a detailed discussion of these computations, which involved several CPU-years of computer time.

Proving results about the distribution of doubly isogenous curves in the \dsixfamily over finite fields seems out of reach.  However, conditional on the Zilber--Pink conjecture we are able to prove that there are finitely many doubly isogenous \dsixcurves in characteristic zero (Theorem~\ref{theorem:finitely many doubly isogenous}).  This involves ruling out the existence of simultaneous isogeny correspondences between elliptic curves appearing in the decomposition of the Jacobians of $C^{[2]}$ and $C^{'[2]}$ (see Section~\ref{sec:computingisogeny}).  We develop an algorithm  to do so, which relies on a criterion of Buium (Theorem~\ref{theorem:buium}) for the existence of isogeny correspondences.  Crucially, this does not require prior knowledge of the degree of the isogeny correspondence.

Finally, following the philosophy of \cite{SutherlandVoloch2019} (extended as in \cite{BooherVoloch2020}), it is natural to attempt to use additional unramified covers to distinguish doubly isogenous curves in our family from one another.  In Sections \ref{sec:covers4} and \ref{sec:triplecovers} we investigate unramified cyclic degree-$4$ covers of \dsixcurves and certain unramified cyclic triple covers, and in Sections \ref{sec:2coverheuristics}--\ref{sec:computingisogeny} we present heuristic and numerical arguments that the isogeny classes of the Jacobians of these covers \emph{can} distinguish non-isomorphic doubly isogenous curves in our family from one another.  Remarkably, these Jacobians also decompose up to isogeny as a product of elliptic curves, which leads to an unexpected computational application, as we now explain.

Call a one-parameter family of abelian varieties defined by equations over $\ZZ[t]$ a ``distinguished family'' if for all sufficiently large primes $p$, values of $t$ that are distinct modulo $p$ yield abelian varieties whose reductions modulo $p$ are non-isogenous.
Poonen~\cite{Poonen2019} has shown that any distinguished family leads to a deterministic polynomial-time algorithm to factor univariate polynomials over finite fields. The existence of a distinguished family remains an open problem, as does the existence of such an algorithm. Sutherland and Voloch \cite{SutherlandVoloch2019} presented a family of genus-$4$ curves with ample numerical evidence suggesting that the Jacobians of these curves form a distinguished family. In Section~\ref{sec:factoring} we present an alternative candidate for a distinguished family, derived from our \dsixfamily, that is much more computationally efficient than the one given in~\cite{SutherlandVoloch2019}. In practice, of course, there are faster probabilistic, polynomial-time algorithms to factor univariate polynomials over a finite field, but none of them are deterministic.

\begin{remark}
This paper and \cite{ArulEtAl2024} study similar questions and share a common history.  An initial version of the work presented here inspired a project at the June, 2020 ICERM Workshop on Arithmetic Geometry, Number Theory, and Computation, which resulted in \cite{ArulEtAl2024}.  That paper considers the one-parameter family of curves of genus $2$ whose automorphism group contains~$D_4$. 
Again, instances of pairs of non-isomorphic doubly isogenous curves were more common than initially expected, even taking into account symmetries of the decomposition of the Jacobian.  In that case, the discrepancies were explained not through the existence of a pair of doubly isogenous curves in characteristic zero but instead through the existence of multiple simultaneous isogeny correspondences between the simple abelian varieties appearing in the decomposition.  The algorithm for finding simultaneous isogeny correspondences given in Section \ref{sec:isogenyalgorithm} can verify that \cite{ArulEtAl2024} found all of them for the $D_4$ family.
\end{remark}

\subsubsection*{Acknowledgements}
Voloch would like to thank S. C. Coutinho for some computational suggestions and J. Pila for some references. We are grateful to the anonymous referee for their careful reading and their detailed comments and suggestions, all of which helped improve the exposition in this paper.

\section{The family of genus-two curves with automorphism group \texorpdfstring{$D_6$}{D6}}
\label{sec:family}

For fields $K$ of characteristic neither $2$ nor~$3$, we will
describe the moduli space for pairs $(C,\eta)$, where $C$ is a genus-$2$ curve
over $K$ and $\eta$ is an embedding of the dihedral group $D_6$ with twelve
elements into $\Aut_K C$, the group of $K$-rational automorphisms of~$C$. Note
that if the automorphism group of a genus-$2$ curve contains $S_3$, then the
subgroup generated by the hyperelliptic involution and $S_3$ is isomorphic
to~$D_6$. 

Such curves have been studied by other authors, including Igusa~\cite{Igusa1960},
Ibukiyama, Katsura, and Oort~\cite{IbukiyamaKatsuraOort1986}, and 
Cardona and Quer~\cite{CardonaQuer2007}. In fact, Cardona and Quer gave a 
classification, over perfect fields~$K$, of genus-$2$ curves whose
\emph{geometric} automorphism groups are isomorphic to~$D_6$; they give a 
geometric moduli space for such curves, describe the field of definition for
each point in the moduli space, and describe the twists of those curves whose
fields of definition are~$K$. For the present work we prefer the parametrization
given below, because we only want to consider curves whose \emph{$K$-rational}
automorphism groups are isomorphic to~$D_6$, and we want to emphasize the
elliptic curves that fall out of our analysis.

We fix a presentation of the group $D_6$:
\begin{equation}
\label{EQ:D12}
D_6 = \langle a, b, i \, | \, a^3 = b^2 = i^2 = (ba)^2 = 1, ia=ai, ib=bi \rangle.
\end{equation}
Note that $i$ is the central involution of $D_6$.

\begin{definition}
\label{D:D6curve}
A \emph{\dsixcurve} over a field $K$ is a pair $(C,\eta)$, where $C$ is a genus-$2$
curve over $K$ all of whose geometric automorphisms are rational over~$K$ and
where $\eta$ is an isomorphism from $D_6$ to the group $\Aut C$ of geometric
automorphisms of~$C$. Two \dsixcurves $(C_1,\eta_1)$ and $(C_2,\eta_2)$ are
\emph{isomorphic} to one another if there is an isomorphism of curves $\gamma\colon C_1\to C_2$
that identifies the isomorphisms $\eta_1$ and~$\eta_2$, that is, such that 
$\eta_2(x) = \gamma\eta_1(x)\gamma^{-1}$ in $\Aut C_2$ for all $x\in D_6$.
\end{definition}

The next lemma shows how to tell whether two $D_6$ structures on the same curve
are isomorphic.

\begin{lemma}
\label{L:tellstructures}
Let $K$ be an arbitrary field and let $(C,\eta)$ and $(C,\eta')$
be two \dsixcurves over $K$ with the same underlying 
genus-$2$ curve $C$. These two \dsixcurves are isomorphic to one another if and
only if $\beta\colonequals \eta(b)$ and $\beta'\colonequals \eta'(b)$ are
conjugate in~$\Aut C$.
\end{lemma}

\begin{proof}
We know that $\beta$ and $\beta'$ are noncentral involutions in 
$\Aut C\cong D_6$. Suppose they are conjugate, with $\gamma\in \Aut C$ an
element with $\gamma\beta\gamma^{-1} = \beta'$. Then $\gamma\beta$ also
conjugates $\beta$ to $\beta'$, and one of the two elements $\gamma$ and
$\gamma\beta$ will conjugate $\eta(a)$ to $\eta'(a)$. Without loss of
generality, we suppose that $\gamma$ does so. Then the isomorphism
$\gamma\colon C\to C$ gives an isomorphism from $(C,\eta)$ to $(C,\eta')$ as
\dsixcurves.

Conversely, if $\gamma\in \Aut C$ is an automorphism that induces an isomorphism
of $(C,\eta)$ to $(C,\eta')$ as \dsixcurves, then 
$\gamma\beta\gamma^{-1} = \beta'$.
\end{proof}

\begin{lemma}
\label{L:twostructures}
Let $(C,\eta)$ be a \dsixcurve over an arbitrary field~$K$.
Up to isomorphism, there is exactly one other
\dsixcurve $(C,\eta')$ with the same underlying~$C$.
\end{lemma}

\begin{proof}
This follows from Lemma~\ref{L:tellstructures} and the fact that the outer 
automorphism group of $D_6$ has order~$2$. In more detail:
If $\alpha\colon D_6 \to D_6$ is a nontrivial outer automorphism, 
then $(C,\eta\alpha)$ is a \dsixcurve that is not isomorphic to $(C,\eta)$, by
Lemma~\ref{L:tellstructures}. Conversely, if $(C,\eta')$ is a \dsixcurve, not
isomorphic to $(C,\eta)$, then $\eta^{-1}\eta'$ is an automorphism of $D_6$
which is necessarily outer, again by Lemma~\ref{L:tellstructures}.
\end{proof}

\begin{notation}
\label{N:basiccurve}
Let $K$ be a field of characteristic not~$2$, and let $r$ and $c$ be elements 
of $K^\times$ with $r\ne \pm 27$ and with
$r\neq 23 \pm 10\sqrt{-2}$. We let $C_{r,c}$ denote the curve over $K$ given by
\begin{equation}
\label{EQ:basiccurve}
cy^2 = x^6 + (r-18)x^4 + (81-2r)x^2 + r.
\end{equation}
Let $\alpha$ and $\beta$ be the automorphisms of $C_{r,c}$ given by
\begin{align*}
\alpha \colon (x,y) &\mapsto \biggl(\frac{x-3}{x+1}, \frac{-8y}{(x+1)^3}\biggr)\\
\beta  \colon (x,y) &\mapsto (-x,y)
\end{align*}
and let $\iota$ be the hyperelliptic involution. We let $\eta_{r,c}$ denote the
map $D_6\to \Aut C_{r,c}$ that takes $a$, $b$, and $i$ to $\alpha$, $\beta$,
and~$\iota$, respectively.
\end{notation}

\begin{remark}
\label{R:toobig}
We exclude the values $r=0$ and $r=27$ in this notation, because for these values
of $r$ the discriminant of the right side of~\eqref{EQ:basiccurve} is~$0$.
We exclude the values $r = -27$ and $r = 23\pm 10\sqrt{-2}$ for a less obvious
reason. Igusa's classification of curves of genus~$2$ with large automorphism
groups~\cite[\S8]{Igusa1960} shows that over an algebraically closed field of
characteristic other than~$2$, there are at most two curves whose automorphism 
groups properly contain~$D_6$. By calculating the Igusa invariants of these curves,
we find that~\eqref{EQ:basiccurve} gives us one of them when~$r=-27$,
and the other when $r = 23\pm 10\sqrt{-2}$. Since we are interested in curves
whose geometric automorphism groups are~$D_6$, we exclude those values of~$r$
in Notation~\ref{N:basiccurve}. For this reason, we introduce the following notation.
\end{remark}

\begin{notation}
Given a field $K$, we denote by $S_K$ the set
$K^\times\setminus\{\pm 27, 23\pm 10\sqrt{-2}\}$.
\end{notation}

Lemma~\ref{L:twostructures} says that a genus-$2$ curve with 
a $D_6$ structure has exactly one other $D_6$ structure that is 
not isomorphic to the first.
Our next result
tells us what the ``partner curve'' of $(C_{r,c},\eta_{r,c})$ is,
if the characteristic of $K$ is neither $2$ nor~$3$.

\begin{lemma}
\label{L:rcpairs}
Let $K$ be a field of characteristic neither $2$ nor $3$, 
let $r$ be an element of $S_K$ and $c$ an element of $K^\times$, and let
$r'\colonequals 729/r$ and $c'\colonequals cr$. The curves $C_{r,c}$ and 
$C_{r',c'}$ are isomorphic to one another, but the \dsixcurves 
$(C_{r,c},\eta_{r,c})$ and  $(C_{r',c'},\eta_{r',c'})$ are not isomorphic to one
another.
\end{lemma}

\begin{proof}
An isomorphism from $C_{r,c}$ to $C_{r',c'}$ is given by
\[
(x,y)\mapsto \Bigl(\frac{3}{x},\frac{27 y}{r x^3}\Bigr).
\]
The pullback of the automorphism $(x,y)\mapsto(-x,y)$ of $C_{r',c'}$ under this
map is the automorphism $(x,y)\mapsto (-x,-y)$ of $C_{r,c}$, which 
is~$\iota\beta$. We check that $\iota\beta$ is not conjugate to $\beta$ in
$\Aut C_{r,c}\cong D_6$, so the two \dsixcurves are not isomorphic to one
another.
\end{proof}

Our main result for this section is that every \dsixcurve over a field of
characteristic neither $2$ nor $3$ is
isomorphic to exactly one of the pairs $(C_{r,c},\eta_{r,c}).$

\begin{proposition}
\label{P:standardform}
Let $K$ be a field of characteristic neither $2$ nor~$3$.
Let $(C,\eta)$ be a \dsixcurve over $K$. Then there is a pair 
$(r,c)\in S_K \times K^\times$, with $r$ unique and $c$ unique up to squares,
such that there is an isomorphism $\theta\colon C\to C_{r,c}$ such that
$\theta^*\eta_{r,c} = \eta$.
\end{proposition}

We present the proof of this proposition later in the section. For now, we begin
with a basic observation about abelian surfaces with an action of~$S_3$. Over
the complex numbers this is an expansion of a special case of a result of 
Ries~\cite[Prop.~4.1]{Ries1997}. While Ries's proof works over more general
fields, for the convenience of the reader we will present a proof rather than
simply cite his work. In the following proposition and its proof, we write endomorphisms of a
product $E_1\times E_2$ of elliptic curves in matrix-like arrays of the form
\[
M = \begin{bmatrix} \alpha & \beta \\ \gamma & \delta\end{bmatrix};
\]
in this notation we mean that 
\begin{align*}
\alpha&\in\End E_1, & \beta&\in\Hom(E_2,E_1),\\
\gamma&\in\Hom(E_1,E_2), & \delta&\in\End E_2,
\end{align*}
and the endomorphism $M$ is given by
$M(P,Q) = (\alpha(P) + \beta(Q), \gamma(P) + \delta(Q)).$

\begin{proposition}
\label{T:S3surfaces}
Let $K$ be a field of arbitrary characteristic, let $A$ be an abelian 
surface over $K$ with a principal polarization~$\lambda$, and let 
$\sigma,\tau\in \Aut_K (A,\lambda)$ be automorphisms such that
\begin{enumerate}
\item $\sigma^3 = \tau^2 = 1$,
\item $\tau\sigma = \sigma^2\tau$, and
\item $\sigma^2 + \sigma + 1 = 0$ in $\End A$.
\end{enumerate}
Then there is a $3$-isogeny $\varphi\colon E_1\to E_2$ of elliptic curves 
over~$K$, unique up to isomorphism, such that there is an isomorphism 
$A\to E_1\times E_2$ that takes the automorphisms $\sigma$ and $\tau$ to the
automorphisms
\begin{equation}
\label{EQ:auts}
\begin{bmatrix} -2 & -\hat{\varphi}\\ \phantom{-}\varphi&\phantom{-}1\end{bmatrix}
\text{\quad and\quad}
\begin{bmatrix} 1 & \phantom{-}\phihat\\ 0&-1\end{bmatrix}
\end{equation}
of $E_1\times E_2$ and that takes the polarization $\lambda$ to the polarization
\begin{equation}
\label{EQ:pol}
\begin{bmatrix}2&\phihat\\ \varphi&2\end{bmatrix}
\end{equation}
of $E_1 \times E_2$.

Conversely, given a $3$-isogeny $\varphi\colon E_1\to E_2$ between two elliptic
curves over~$K$, the automorphisms $\sigma$ and $\tau$ of the surface
$E_1\times E_2$ given by the two matrices~\eqref{EQ:auts} satisfy conditions 
$(1)$, $(2)$, and~$(3)$, and the automorphisms respect the principal
polarization of $E_1\times E_2$ given by~\eqref{EQ:pol}.
\end{proposition}

\begin{proof}
Let $R$ be the subring $\ZZ[\sigma,\tau]$ of $\End A$. We define several
elements of~$R$:
\begin{align*}
\eps_1 &\colonequals - (1 + \tau)\sigma^2\\
\eps_2 &\colonequals -\sigma (1 - \tau)\\
\Phi &\colonequals (1 - \sigma^2)\eps_1\\
\Phi' &\colonequals (1 - \sigma)\eps_2.
\end{align*}
All of these elements act on $A$. We calculate that $\eps_1$ and $\eps_2$ are
orthogonal idempotents that sum to~$1$. Note that $\tau\ne\pm 1$ because $\tau$
does not commute with $\sigma$, so $\eps_1$ and $\eps_2$ are nontrivial
idempotents. For each $i$ we let $E_i$ be the elliptic curve $\eps_i A$.

We check that $\eps_1\Phi\eps_1 = \eps_2\Phi\eps_2 = \eps_1\Phi\eps_2 = 0$.
Let $\varphi\colonequals \eps_2\Phi\eps_1\in\Hom(E_1,E_2)$, so that under the
isomorphism $A \cong E_1\times E_2$ provided by the idempotents, we have
\[
\Phi = \begin{bmatrix} 0 & 0\\ \varphi & 0\end{bmatrix}.
\]
Similarly, if we let $\varphi'\colonequals \eps_1\Phi'\eps_2\in\Hom(E_2,E_1)$,
then
\[
\Phi' = \begin{bmatrix} 0 & \varphi' \\ 0 & 0\end{bmatrix}.
\]
Computing the values of $\eps_i \sigma\eps_j$ and  $\eps_i \tau\eps_j$, we find
that under the isomorphism $A\cong E_1\times E_2$ we have
\[
\sigma =  \begin{bmatrix} -2 & -\varphi'\\ \phantom{-}\varphi&\phantom{-}1\end{bmatrix}
\text{\quad and\quad}
\tau =  \begin{bmatrix} 1 & \phantom{-}\varphi'\\ 0&-1\end{bmatrix},
\]
which will be \eqref{EQ:auts} if we can show that $\varphi'=\phihat$.
We check that $\eps_1\varphi'\varphi\eps_1 = 3\eps_1$ and
$\eps_2\varphi\varphi'\eps_2 = 3\eps_2$, so \emph{either} $\varphi$ is an
isogeny of degree~$3$ and $\varphi'$ is its dual, \emph{or} $E_1\cong E_2$ and
we can choose the isomorphism so that $\varphi = 1$ and $\varphi' = 3$ or vice
versa.

Now we use the assumption (unused until this point) that $A$ has a principal
polarization~$\lambda$. Under the isomorphism $A\cong E_1\times E_2$ we must
have
\[
\lambda =  \begin{bmatrix} a & \psihat\, \\ \psi &b\end{bmatrix} 
\]
for some morphism $\psi\colon E_1\to E_2$ with dual $\psihat$ and for some
positive integers $a$ and $b$. The assumption that $\lambda$ is principal means
that $ab - \psihat\psi = 1$. For the automorphisms $\sigma$ and $\tau$ to
respect the polarization, we need 
\begin{align}
\label{EQ:sigmarespect} \lambda &= \hat{\sigma} \lambda \sigma\\
\label{EQ:taurespect}  \lambda &= \hat{\tau} \lambda \tau,
\end{align}
where $\hat\sigma$ and $\hat{\tau}$ are the dual isogenies of $\sigma$
and~$\tau$. Note that in the matrix representation for elements of 
$\End(E_1\times E_2)$, the dual isogeny is just the conjugate transpose.

Suppose we are in the case where $E_1 = E_2$ and $\varphi = 1$ and
$\varphi' = 3$. Then the condition that $\lambda = \hat{\tau} \lambda \tau$ can
be written
\[
\begin{bmatrix} a & \psihat\, \\ \psi &b\end{bmatrix}
=
\begin{bmatrix} 1 & \phantom{-}0\\ 3&-1\end{bmatrix}
\begin{bmatrix} a & \psihat\, \\ \psi &b\end{bmatrix}
\begin{bmatrix} 1 & \phantom{-}3\\ 0&-1\end{bmatrix},
\]
which simplifies to
\[
\begin{bmatrix} a & \psihat\, \\ \psi &b\end{bmatrix}
=
\begin{bmatrix}
a& 3a - \psihat\\
3a - \psi  & 9a - 3\psihat-3\psi + b
\end{bmatrix}.
\]
From this equality we deduce that $2\psi = 2\psihat = 3a.$ Similarly, from the 
condition $\lambda = \hat{\sigma}\lambda\sigma$ we deduce that we also have
$b = 3a$. Combining these equalities with the condition that
$ab-\psihat\psi = 1$ we find that $9a^2 - 9a^2/4 = 1$, so $a^2 = 4/27$, an 
impossibility.

If we assume that $\varphi=3$ and $\varphi'=1$ we again reach a contradiction.
Therefore, we must be in the case where $\varphi$ is a $3$-isogeny and where 
$\varphi'$ is the dual isogeny $\phihat$ of $\varphi$. In this case, the 
equalities $\lambda = \hat{\sigma}\lambda\sigma$ and 
$\lambda = \hat{\tau} \lambda \tau$ and $ab - \psihat\psi = 1$ determine 
$\lambda$ completely: We find that
\[
\lambda = \begin{bmatrix}2&\phihat\\ \varphi&2\end{bmatrix}.
\]
This proves the first part of the theorem.

To prove the second part of the theorem we simply observe that if we take 
$\sigma$, $\tau$, and $\lambda$ as in \eqref{EQ:auts} and \eqref{EQ:pol}, then 
\eqref{EQ:sigmarespect} and \eqref{EQ:taurespect} hold, and $\sigma$ and $\tau$
satisfy $\sigma^3 = \tau^2 = 1$, $\tau\sigma = \sigma^2\tau$, and 
$1 + \sigma + \sigma^2 = 0$.
\end{proof}

The final ingredient for our proof of Proposition~\ref{P:standardform} can be
interpreted as an explicit description of the modular curve $X_0(3)$ over fields
of characteristic neither $2$ nor~$3$.

\begin{lemma}
\label{L:X0(3)}
Let $K$ be a field of characteristic neither $2$ nor~$3$.
Let $E$ be an elliptic curve over $K$ and let $G$ be a $K$-rational finite 
subgroup scheme of $E$ of {order}~$3$. If $j(E)$ and $j(E/G)$ are not both~$0$, then
there is a pair $(r,c)\in K^\times\times K^\times$, with $r$ unique and $c$
unique up to squares, such that there is an isomorphism from $E$ to the curve
\begin{equation}
\label{EQ:basicec}
c y^2 =  x^3 + (r-18)x^2 + (81-2r)x + r
\end{equation}
that takes $G$ to the subgroup scheme of {order}~$3$ whose nonzero points
satisfy $x=-3$. For this pair $(r,c)$, the quotient curve $E/G$ is isomorphic to
\begin{equation}
\label{EQ:ecquotient}
c y^2 = x^3 + (81 - 2r)x^2 + r(r - 18)x + r^2.
\end{equation}

Conversely, given a pair $(r,c)\in K^\times\times K^\times$ with $r\ne 27$,
equation~\eqref{EQ:basicec} defines an elliptic curve $E$ such that the points
on $E$ with $x=-3$ are the nonzero elements of a subgroup scheme $G$ of order~$3$, and the
$j$-invariants of $E$ and $E/G$ are not both~$0$.
\end{lemma}

\begin{proof}
First let us consider the case where $j(E)\ne 0$. We can write $E$ in the form
$y^2 = x^3 + ax^2 + bx + d$, and by shifting coordinates we may assume that the
(shared) $x$-coordinate of the nonzero points of $G$ is~$0$. Since the
$3$-division polynomial of this curve is 
\[
3x^4 + 4ax^3 + 6bx^2 + 12dx + 4ad - b^2,
\]
this means that $b^2 = 4ad$. We check that then $a$ and $b$ are both nonzero,
since otherwise the curve would either be singular or have $j$-invariant~$0$.
Let $c = -8a/b$, and replace $x$ with $x/c$ and $y$ with $y/c$. Then we find
that $cy^2 = x^3 + a' x^2 + b'x + d'$, where $a' = ac$ and $b'=bc^2 = -8a'$ and 
$d' = dc^3 = 16 a'$. If we set $r = a'+27$ and replace $x$ with $x+3$, we find
that $E$ is isomorphic to the curve~\eqref{EQ:basicec}, and the $x$-coordinate
of the nonzero points in $G$ is $x = -3$. This normalization process shows that
$r$ is unique, and $c$ is unique up to squares.

If $j(E) = 0$ then $E$ can be written in the form $y^2 = x^3 + d$. The 
$3$-division polynomial of this curve is $3x(x^3+4d)$, and the $3$-torsion
points whose $x$ coordinates satisfy $x^3 + 4d$ give rise to quotient curves
with nonzero $j$-invariants. If we let $a$ be the $x$-coordinate of the nonzero
points in $G$ then $d = -a^3/4$. Let $c = -8/a$, and replace $x$ with $x/c$ and
$y$ with $y/c$. Then we find that $cy^2 = x^3 + 128$. If we replace $x$ with
$x - 5$ and set $r = 3$, we find that $E$ is isomorphic to the 
curve~\eqref{EQ:basicec}, and the $x$-coordinate of the nonzero points in $G$ is
$x = -3$. Once again, the normalization process shows that $r$ is unique, and
$c$ is unique up to squares.

A $3$-isogeny from~\eqref{EQ:basicec} to~\eqref{EQ:ecquotient} with
kernel $G$ is given by
\begin{align*}
(x,y) &\mapsto \Bigl( \frac{x^3 + (r - 27)x^2 + (-10r + 243)x + (25r - 729)}{(x + 3)^2},\\
      &\qquad\qquad      \frac{x^3 + 9x^2 + (16r - 405)x + (-80r + 2187)}{(x+3)^3} \, y\Bigr).
\end{align*}

The converse statement is easy to check; the restriction $r\ne 27$ is necessary
in order for~\eqref{EQ:basicec} and \eqref{EQ:ecquotient} to define nonsingular 
curves.
\end{proof}

\begin{remark}
\label{R:3isogeny}
Equation~\eqref{EQ:basicec} can be rewritten as
$y^2 = (x+3)^3 + (r-27)(x-1)^2.$ More generally, one can check that if an 
elliptic curve over $K$ in Weierstrass form $y^2=f$ has a rational $3$-isogeny
such that the $x$-coordinate of the nonzero points in the kernel is~$a$, then
the polynomial $f$ can be written $f = (x - a)^3 + b(x - c)^2$ for some $b,c\in K$.
Conversely, every curve of this form has a rational $3$-isogeny.
\end{remark}

\begin{proof}[Proof of Proposition~\textup{\ref{P:standardform}}]
Let $(J,\lambda)$ be the canonically principally polarized Jacobian of $C$ and let $\sigma$
and $\tau$ be the automorphisms $\eta(a)^*$ and $\eta(b)^*$ of~$(J,\lambda)$. 
Applying Theorem~\ref{T:S3surfaces}, we find that there is a $3$-isogeny 
$\varphi\colon E_1\to E_2$ of elliptic curves over $K$ and an isomorphism 
$J\cong E_1\times E_2$ that identifies $\sigma$ and $\tau$ with the 
automorphisms in~\eqref{EQ:auts} and that identifies $\lambda$ with the 
polarization given by~\eqref{EQ:pol}.

Consider the degree-$4$ endomorphism $\Phi$ of $E_1\times E_2$ given by
\[
\Phi\colonequals
\begin{bmatrix}
\phantom{-}1  & -\phihat \\
-\varphi      & \phantom{-}1
\end{bmatrix}.
\]
We note that $\Phi$ is equal to its own dual isogeny (when we identify
$E_1\times E_2$ with its dual via the canonical isomorphism of each factor with
its dual), and we have 
\[
\begin{bmatrix}
2 & 0  \\
0 & 2
\end{bmatrix}
=
\hat\Phi
\begin{bmatrix}
2       & \phihat\\ 
\varphi &  2
\end{bmatrix}
\Phi.
\]
We check that the kernel of $\Phi$ is the graph of the isomorphism 
$\psi\colon E_1[2]\to E_2[2]$ induced by the isogeny~$\varphi$. 
By~\cite[Proposition~3, p.~324]{HoweLeprevostPoonen2000}, we see that $C$ is
therefore the curve obtained by ``gluing'' $E_1$ to $E_2$ along their 
$2$-torsion subgroups via~$\psi$. In particular, this implies that $\psi$ must
not be the graph of an isomorphism $E_1\to E_2$ restricted to $E_1[2]$, and a
short computation shows that this means precisely that either $j(E_1)\ne 0$ or
$j(E_2)\ne 0$.

Applying Lemma~\ref{L:X0(3)}, we find that there are elements $r$ and $c$ of $K$
such that $E_1$ is isomorphic to $cy^2 =  x^3 + (r-18)x^3 + (81-2r)x + r$ and 
$E_2$ is isomorphic to $c y^2 = x^3 + (81-2r)x^2 + r(r - 18)x + r^2$, and by 
choosing the isomorphisms properly we can assume that the isogeny
$\varphi\colon E_1\to E_2$ is given by the formula in the proof of 
Lemma~\ref{L:X0(3)}. We check that this isogeny induces the isomorphism 
$E_1[2]\to E_2[2]$ given by sending a nontrivial element $(s,0)\in E_1[2]$ to 
the element $(r/s,0)$ of $E_2[2]$. Using the formulas 
from~\cite[Proposition~3, p.~324]{HoweLeprevostPoonen2000} and simplifying, we
find that gluing $E_1$ and $E_2$ together in this way gives us the curve
$C_{r,c}$, so we have found an isomorphism $\theta\colon C \to C_{r,c}$.

By Lemmas~\ref{L:twostructures} and~\ref{L:rcpairs}, it follows that
the \dsixcurve $(C,\eta)$ is isomorphic either to $(C_{r,c}, \eta_{r,c})$ or to 
$(C_{r',c'}, \eta_{r',c'})$, where $r'=729/r$ and $c'=cr$.
\end{proof}

In the course of proving Proposition~\ref{P:standardform} we computed the isogeny
decomposition of the Jacobian of $C_{r,c}$ (compare to~\cite[Remark~1.4, p.~132]{IbukiyamaKatsuraOort1986}):

\begin{corollarynote}
\label{P:baseE}
Let $K$ be a field of characteristic neither $2$ nor~$3$.
If $r\in S_K$ and $c\in K^\times$, the Jacobian of $C_{r,c}$ is isogenous to the square of the elliptic curve
\[
E_{r,c}\colon\quad c y^2 = x^3 + (r-18)x^2 + (81-2r)x + r,
\]
whose $j$-invariant is $j_r\colonequals -(r-3)^3 (r-27)/r$.
\end{corollarynote}

Proposition~\ref{P:standardform} suggests the following definition.

\begin{definition}
\label{D:invariant}
Let $C$ be a \dsixcurve over a field of characteristic neither $2$ nor~$3$.
The \emph{coarse invariant} of $C$ is the value $r + 729/r$,
where $r$ is the element of $S_K$ such that there is a $c$ with $C \cong C_{r,c}$
as \dsixcurves.
\end{definition}

Note that if $K$ is algebraically closed, two \dsixcurves over $K$ are isomorphic
to one another as curves if and only if their coarse invariants are equal.

\section{The \texorpdfstring{$2$}{2}-torsion and unramified elementary abelian \texorpdfstring{$2$}{2}-covers}
\label{sec:covers2}

In this section we will look at the subfamily of \dsixcurves $C_{r,c}$ that have all of their Weierstrass points rational, and we will
write down explicit models of the Prym varieties of their unramified double
covers. Our first result, whose proof we leave to the reader, tells us which
values of $r$
lead to curves with rational Weierstrass points. Throughout this section, $K$
will be a field of characteristic neither $2$ nor~$3$.

\begin{lemma}
\label{L:rationalW}
Suppose the curve $C_{r,c}$ has a rational Weierstrass point, say $P = (u,0)$.
Then 
\[
r = - \frac{u^2 (u^2 - 9)^2}{(u^2 - 1)^2},
\]
and all of the Weierstrass points of $C_{r,c}$ are rational. In this case, the 
$x$-coordinates of the rational Weierstrass points are
\[
\pm u,  \qquad \pm \frac{u + 3}{u - 1}, \quad \text{and\quad}\pm \frac{u - 3}{u + 1}.
\]
The automorphism group of $C_{r,c}$ acts transitively on these six points.\qed
\end{lemma}

\begin{notation}
\label{N:ucurves}
Let $U_\bad\subseteq K$ be the set $U_\bad \colonequals U_\sing \cup U_\toobig$, 
where
\begin{align*}
U_\sing &\colonequals 
K \cap \bigl\{
0, \ 
\pm 1, \  
\pm 3, \ 
\pm\sqrt{-3}
\bigr\} \qquad\textup{and}\\
U_\toobig &\colonequals  
K \cap \bigl\{
\pm\sqrt{3}, \ 
\pm3 \pm 2\sqrt{3}, \ 
\pm \sqrt{-1} \pm \sqrt{2}, \ 
\pm (1 \pm \sqrt{-2})(1 \pm \sqrt{2}) \bigr\}.
\end{align*}
Given $u\in K\setminus U_\bad$ and $c\in K^\times$, set
$r \colonequals -u^2 (u^2 - 9)^2/(u^2 - 1)^2.$ We let $D_{u,c}$ denote the curve
$C_{r,c}$ from Notation~\textup{\ref{N:basiccurve}}.
\end{notation}

The elements of $U_\sing$ are the values of $u$ that produce singular curves;
in particular, $u=\pm 1$ gives $r=\infty$; 
$u=0$ and $u=\pm 3$ give $r=0$, which leads to a singular curve; and
$u=\pm\sqrt{-3}$ gives $r=27$, which also leads to a singular curve.
The elements of $U_\toobig$ lead to $r=-27$ or $r = 23\pm 10\sqrt{-2}$, which 
give curves with automorphism groups larger than~$D_6$.

A computation that we leave to the reader shows that there are twelve pairs 
$(u,c)$ that give rise to a given curve $D_{u,c}$:

\begin{lemma} \label{lemma:isomorphic D_{u,c}}
If $(u,c)$ and $(u',c')$ are elements of 
$(K\setminus U_\bad) \times (K^\times/K^{\times 2})$ such that 
$D_{u',c'}$ is isomorphic as a curve to $D_{u,c}$, then $(u',c')$ is one of the
following, where we understand the second element to be given only up to
squares in $K^\times$\textup{:}
\begin{align*}
&(\pm u, c),&  
&\qquad \Bigl( \frac{\pm (u + 3)}{u - 1},c\Bigr), &
&\qquad \Bigl( \frac{\pm (u - 3)}{u + 1},c\Bigr),\\
&\Bigl( \frac{\pm 3}{u},-c\Bigr), &
&\qquad \Bigl( \frac{\pm 3(u-1)}{u + 3},-c\Bigr), &
&\qquad \Bigl( \frac{\pm 3(u+1)}{u - 3},-c\Bigr).
\end{align*}
If $u$ corresponds to a given value $r$ as in Notation~\textup{\ref{N:ucurves}},
the pairs on the first line correspond to this same value of $r$, while the
pairs on the second line correspond to the value $729/r$.\qed
\end{lemma}

\begin{notation}
\label{N:Wpoints}
Given $u \in K\setminus U_\bad$ and $c\in K^\times$, let us write the
$x$-coordinates of the Weierstrass points of $D_{u,c}$ as 
\begin{align*}
u_1 &\colonequals \phantom{-}u & u_3 &\colonequals \phantom{-}(u + 3)/(u - 1) & u_5 &\colonequals \phantom{-}(u - 3)/(u + 1)\\
u_2 &\colonequals -u & u_4 &\colonequals -(u + 3)/(u - 1) & u_6 &\colonequals -(u - 3)/(u + 1),
\end{align*}
so that an equation for $D_{u,c}$ is
\[
cy^2 = (x - u_1)(x - u_2)(x - u_3)(x - u_4)(x - u_5)(x - u_6).
\]
\end{notation}

The unramified double covers of $D_{u,c}$ are the curves obtained by adjoining
to the function field of $D_{u,c}$ a root of an equation of the form
\begin{equation}
\label{EQ:cover}
d z^2 = (x - u_i)(x - u_j),
\end{equation}
where $1\le i < j\le 6$ and where $d$ is nonzero. (This follows from class field
theory for function fields \cite[Chapter VI]{Serre1988}; the details for the 
analogous case of hyperelliptic Riemann surfaces are spelled out 
in~\cite{FuertesGonzalez-Diez2006}.) We can normalize (up to squares)
the values of $d$ we consider by demanding that a fixed Weierstrass point 
of $D_{u,c}$ split in the cover; when we consider the set of all unramified
double covers, it does not matter which Weierstrass point we pick to normalize
our covers, because the automorphism group of $D_{u,c}$ acts transitively on 
the Weierstrass points. In what follows, we will demand that the point
$P_1 \colonequals (u,0)$ split in our extensions. 

\begin{notation}
\label{N:doublecovers}
We let $\Dbar_{u,c,i,j}$ denote the curve defined by \eqref{EQ:cover}, normalized
as described above, and we let $\zeta_{u,c,i,j}$ denote the degree-$2$ morphism
$\Dbar_{u,c,i,j}\to D_{u,c}$. If the pair $(u,c)$ is clear from context, we
also write $\zeta_{i,j}\colon \Dbar_{i,j}\to D$ for this cover.
\end{notation}

The group $\Aut D_{u,c}$ acts on the Weierstrass points of $D_{u,c}$. We
use this to define an action of $\Aut D_{u,c}$ on the integers $i$ between $1$ 
and $6$ by demanding that for $\gamma \in \Aut D_{u,c}$ we have
$u_{\gamma(i)} = \gamma(u_i)$; here we use the numbering of the
Weierstrass points given in Notation~\ref{N:Wpoints}. This gives us an action of
$\Aut D_{u,c}$ on the double covers $\zeta_{i,j}$, defined by 
$\gamma(\zeta_{i,j}) = \zeta_{\gamma(i),\gamma(j)}$. We check that 
$\gamma(\zeta_{i,j})$ is geometrically isomorphic to the double cover
$\zeta_{i,j}\circ \gamma^{-1}$. The two covers are not necessarily isomorphic
over~$K$, because in the second cover the point $P_1$ may not split.

The Prym variety $F_{i,j}$ of $\zeta_{i,j}$ is the Jacobian of the genus-$1$ curve
\[
cd y^2 = \prod_{k\ne i,j} (x - u_k).
\]
Each $F_{i,j}$ can be written in the form $e y^2 = x(x-1)(x-\lambda)$ for
certain values of $e$ and~$\lambda$. A straightforward computation gives the
information in Table~\ref{table:lambdas}, which specifies the values of
$\lambda$ and $e$ for each Prym $F_{i,j}$. These Pryms fall into
orbits under the action of $\Aut D_{u,c}\cong D_6$, and in 
Table~\ref{table:lambdas} we group the data by orbit. There is one orbit of 
size~$6$, and there are three of size~$3$.

\begin{table}[t]
\caption{Data for the Prym varieties $F_{i,j}$ for the unramified double covers 
$\zeta_{u,c,i,j}$ of $D_{u,c}$, grouped into orbits under the action of~$D_{6}$.
The first column specifies the name we give to the orbit, the second
gives the indices of the Weierstrass points that specify the $2$-torsion point
on the Jacobian corresponding to the cover, and the third and fourth columns 
give values of $\lambda$ and $e$ such that the Prym is isomorphic to
$e y^2 = x(x-1)(x-\lambda)$.}
\label{table:lambdas}
\begin{tabular}{c c c l}
\toprule
Orbit & & &\\
label & $(i,j)$ & $\lambda$ & \qquad$e$\\
\midrule
\multirow{6}{1em}{$6$}  & $(1, 4)$ & \multirow{6}{9em}{\hfill$\displaystyle\frac{4u}{(u - 1)(u + 3)}$\hfill}
                                                       & $-u(u - 3)(u + 1)$\\
  & $(1, 5)$ &                                         & $u(u - 3)(u + 1)$\\
  & $(2, 3)$ &                                         & $-cu(u - 3)(u - 1)(u + 1)(u + 3)(u^2 + 3)$\\
  & $(2, 6)$ &                                         & $cu(u^2 + 3)$\\
  & $(3, 6)$ &                                         & $c(u - 3)(u + 1)(u^2 + 3)$\\
  & $(4, 5)$ &                                         & $c(u - 1)(u + 3)(u^2 + 3)$\\[2ex] 

  & $(1, 2)$ & \multirow{3}{9em}{\hfill$\displaystyle\frac{16 u^2}{(u^2 + 3)^2}$\hfill}
                                                       & $-(u - 1)(u + 1)(u + 3)(u - 3)$\\
3a& $(3, 5)$ &                                         & $c(u - 3)(u + 1)(u^2 + 3)$\\
  & $(4, 6)$ &                                         & $c(u - 1)(u + 3)(u^2 + 3)$\\[2ex]

  & $(1, 3)$ & \multirow{3}{9em}{\hfill$\displaystyle\frac{(u - 3)^2(u + 1)^2}{(u^2 + 3)^2}$\hfill} 
                                                       & $-u(u - 1)(u + 3)$\\
3b& $(2, 4)$ &                                         & $-cu(u^2 + 3)$\\
  & $(5, 6)$ &                                         & $-c(u - 1)(u + 3)(u^2 + 3)$\\[2ex]

  & $(1, 6)$ & \multirow{3}{9em}{\hfill$\displaystyle\frac{(u - 1)^2(u + 3)^2}{(u^2 + 3)^2}$\hfill} 
                                                       & $u(u - 3)(u + 1)$\\
3c& $(2, 5)$ &                                         & $cu(u^2 + 3)$\\
  & $(3, 4)$ &                                         & $-c(u - 3)(u + 1)(u^2 + 3)$\\
\bottomrule

\end{tabular}
\end{table}

We also calculate the $j$-invariants of the elliptic curves in each orbit:
\begin{align*}
j_6 &= \frac{16 (u^2 + 3)^6}{u^2 (u-3)^2 (u-1)^2 (u + 1)^2 (u + 3)^2} \\[1ex]
j_{\threea} &= \frac{(u^8 - 4 u^6 + 214 u^4 - 36 u^2 + 81)^3}
                 { u^4   (u^2 + 3)^4   (u - 3)^2   (u - 1)^2   (u + 1)^2   (u + 3)^2 } \\[1ex]
j_{\threeb} &= \frac{16 (u^8 - 4 u^7 + 20 u^6 + 52 u^5 - 26 u^4 - 156 u^3 + 180 u^2 + 108 u + 81)^3}
                 {u^2   (u^2 + 3)^4   (u - 3)^4   (u - 1)^2   (u + 1)^4   (u + 3)^2 }\\[1ex]
j_{\threec} &= \frac{16 (u^8 + 4 u^7 + 20 u^6 - 52 u^5 - 26 u^4 + 156 u^3 + 180 u^2 - 108 u + 81)^3}
                 {u^2   (u^2 + 3)^4   (u - 3)^2   (u - 1)^4   (u + 1)^2   (u + 3)^4 }.
\end{align*}
We note that $j_6$ can also be expressed in terms of $r$; we have
$j_6 = -16 (r-27)^2 / r$.

Let $D^{[2]}_{u,c}\to D_{u,c}$ be the maximal abelian $2$-group extension of $D_{u,c}$ 
in which $P_1$ splits completely. This extension is the compositum of the fifteen double
covers given above, and applying a decomposition theorem of Kani 
and Rosen~\cite[Theorem~C]{KaniRosen1989}, we obtain the following proposition.

\begin{proposition}
\label{P:2decomp}
We have
\[
\Jac D^{[2]}_{u,c} \sim
\Jac D_{u,c} \times \prod F_{i,j},
\]
where the product is over all unordered pairs $(i,j)$ of distinct integers 
between $1$ and~$6$.

Let $F_6$, $F_\textup{\threea}$, $F_\textup{\threeb}$, and $F_\textup{\threec}$ be elliptic curves with
$j$-invariants $j_6$, $j_\textup{\threea}$, $j_\textup{\threeb}$, and $j_\textup{\threec}$, respectively,
and let $E_{r,c}$ be the elliptic curve from Corollary~\textup{\ref{P:baseE}}.
Then over the algebraic closure of $K$ we have
\[
\Jac D^{[2]}_{u,c} \sim E_{r,c}^2 \times 
                   F_6^6 \times F_\textup{\threea}^3 \times F_\textup{\threeb}^3 \times F_\textup{\threec}^3,
\]
so that only five elliptic curves appear in the geometric isogeny decomposition 
of $\Jac D^{[2]}_{u,c}$.\qed
\end{proposition}

\section{Unramified cyclic \texorpdfstring{$4$}{4}-covers}
\label{sec:covers4}

In this section we sketch a method for computing the zeta functions of all 
unramified cyclic degree-$4$ extensions of \dsixcurves. 

Let $C$ be a \dsixcurve over an algebraically closed field $K$ of
characteristic not~$2$ or~$3$, so that $C\cong D_{u,c}$ as in Notation~\ref{N:ucurves}
for some $u,c\in K$. The $4$-torsion subgroup of the Jacobian $J$ of $C$ is isomorphic
to $(\ZZ/4\ZZ)^4$, so there are $120$ cyclic order-$4$ 
subgroups of~$J$. By geometric class field theory, there are $120$ 
unramified degree-$4$ covers of $C$ corresponding to these cyclic subgroups.
Every such degree-$4$ cover $D_2 \to C$ factors uniquely as a composition of
unramified degree-$2$ covers $D_2 \to D_1 \to C$.  
The cover $D_1 \to C$ corresponds to the unique point of order $2$ in
the cyclic order-$4$ subgroup of $J$, and is represented by a pair of
Weierstrass points of $C$.  We may choose a model for $C$ in which these points
lie above the points $0$ and $\infty$ of the $x$-line, so that we
can write $C$ in the form
\[
y^2 = x(x-a_1)(x-a_2)(x-a_3)(x-a_4)
\]
for four distinct nonzero elements $a_i$ of~$K$. 
The unramified double cover $D_1$ is obtained by setting $x = z^2$ and $w = y/z$, 
so we have
\begin{equation}
\label{EQ:halfway}
w^2 = (z^2-a_1)(z^2-a_2)(z^2-a_3)(z^2-a_4).
\end{equation}
Choose square roots $b_i$ of the $a_i$, and consider the unramified double
cover $D_2$ of the curve $D_1$ obtained by adjoining a square root of
$(z-b_1)\cdots(z-b_4)$, say
\begin{equation}
v^2 = (z-b_1)(z-b_2)(z-b_3)(z-b_4).
\end{equation}
We observe that the composite unramified cover $D_2\to D_1\to C$ is cyclic: 
One order-$4$ automorphism of $D_2\to C$ is given by
\[
(z,w,v)\mapsto(-z, -w,  w/v).
\]
The cover $D_2\to D_1$ arising from one choice of square roots $\{b_i\}$ is 
isomorphic to the one obtained by negating all the square roots, but otherwise 
different choices of square roots give different covers. Thus, this construction
gives $8$ different cyclic degree-$4$ extensions of $C$ that factor through the
degree-$2$ cover $D_1\to C$. But there are exactly $8$ cyclic order-$4$
subgroups of $J(K)$ that contain a given $2$-torsion point, so we have constructed all such covers.

This construction also shows how the $2$-dimensional Prym 
variety of $D_2\to D_1$ decomposes up to isogeny. The involutions $j_1$ and
$j_2$ of $D_2$ given by $(z,w,v)\mapsto (z,-w,v)$ and $(z,w,v)\mapsto (z,-w,-v)$ 
commute with one another, and the quotient of $D_2$ by the group generated by
$j_1$ and $j_2$ has genus~$0$. The quotient of $D_2$ by the involution $j_1j_2$
is $D_1$; the quotients of $D_2$ by the involutions $j_1$ and $j_2$ are the 
genus-$1$ curves
\begin{equation}
\label{EQ:quarticPrymfactor}
v^2 = (z-b_1)(z-b_2)(z-b_3)(z-b_4)
\end{equation}
and
\[
u^2 = (z+b_1)(z+b_2)(z+b_3)(z+b_4),
\]
respectively, where $u = w/v$. From~\cite[Theorem~C]{KaniRosen1989} we see that
the Jacobian of $D_2$ is isogenous to the product of the Jacobian of $D_1$ and
the Jacobians of these two genus-$1$ curves, and hence the Prym of $D_2\to D_1$ is
isogenous to the product of the Jacobians of the genus-$1$ curves. These two curves
are isomorphic to one another, so in fact the Prym is isogenous to the square of
the Jacobian of the genus-$1$ curve given by~\eqref{EQ:quarticPrymfactor}.

Now suppose our base field $K$ is a finite field. Once again, assume that $C$ is
a \dsixcurve over $K$ with all of its Weierstrass points rational. Then we have
$C\cong D_{u,c}$ for some $u,c\in K$; we fix one such $u$ and $c$ now.
We use the Weierstrass point $W=(u,0)$ to be the base point defining an 
Abel--Jacobi map from $C$ to $J$; all of the unramified covers of $C$ that
we consider will be chosen so that the point $W$ splits completely.

Let $C^{[4]} \to C$ be the degree-$256$ extension obtained from pulling back the
multiplication-by-$4$ map on $J$ via the Abel--Jacobi embedding of $C$ into~$J$.
Over the algebraic closure of~$K$, this extension is the compositum of all of 
the unramified degree-$4$ extensions of $C$ constructed above. The construction 
makes it clear that each of these extensions can be defined over the quadratic 
extension $L$ of~$K$. Since the compositum of these degree-$4$ covers is
$K$-rational, every non-rational elliptic factor in the Prym of one degree-$4$ 
cover has its Galois conjugate appear as a factor in the Prym of another cover.

The cover $C^{[4]} \to C$ factors through the cover $C^{[2]} \to C$ defined
just before Proposition~\ref{P:2decomp}. Since  $C^{[4]}$ has genus $257$ and
$C^{[2]}$ has genus~$17$ (by the Riemann--Hurwitz formula), the Prym variety $P_{4,2}$
of the cover $C^{[4]} \to C^{[2]}$ has dimension~$240$. Using 
\cite[Theorem~C]{KaniRosen1989} once again, we find that this Prym 
decomposes up to isogeny as the product of the Prym varieties of the double
covers $D_2\to D_1$ obtained from the unramified degree-$4$ extensions $D_2\to C$
constructed above. Therefore, $P_{4,2}$ decomposes up to isogeny as the 
product of $240$ elliptic curves.

\begin{lemma}
\label{L:4action}
The automorphism group of $C$ acts on the set of unramified cyclic degree-$4$
covers of $C$ over the algebraic closure of~$K$, and the $240$ elliptic curves 
described above fall into $26$ orbits under this action. 
\end{lemma}

\begin{proof}
Class field theory for curves says that the unramified cyclic degree-$4$ covers
correspond to the index-$4$ subgroups of $(\Jac C)[4]$ with cyclic quotients.
Since the automorphism group of $C$ acts on the set of such subgroups, it
acts on the set of covers as well. Using the explicit description of these
covers given above, one can calculate the orbits of this action. We leave it to
the reader to carry out this calculation and to verify that there are $26$ orbits.
\end{proof}

These $26$ orbits provide us with more data associated to any \dsixcurve $C$
with rational Weierstrass points, say $C \cong D_{u,c}$. By keeping track of
rationality conditions, we can produce a list of the $240$ elliptic curves over 
$L$ attached to~$C$, and of their orbits under the action of the automorphism
group. To each such elliptic curve $E$ we can associate a $\lambda$-invariant 
$\lambda\in L$ and a twisting factor $e\in L$ such that $E$ is given by 
$y^2 = ex(x-1)(x-\lambda)$, and we can also associate a set of expressions in 
the parameters $c$ and $u$ that must be squares in~$K$ in order for the elliptic
factor to be $K$-rational. Such a list is available in the file 
\texttt{fourfold.magma} in the GitHub repository associated to this paper
\cite{DoublyIsogenousRepo}.

Using this list we can compute even more data to associate with $C\cong D_{u,c}$:
For any given values of $u$ and $c$ it is an easy matter to compute the Weil
polynomial of the Prym variety of the cover $C^{[4]}\to C^{[2]}$.
For the $K$-rational elliptic factors, we obtain a factor consisting of the Weil 
polynomial of the elliptic curve. The elliptic factors defined only over the 
quadratic extension of $K$ come in pairs, and for each pair we get a factor
consisting of the Weil polynomial of the restriction of scalars of one of the
curves. For more details, see the discussion of the analogous situation 
in Section~\ref{subsec:general triple}.

\section{Unramified cyclic \texorpdfstring{$3$}{3}-covers}
\label{sec:triplecovers}

In this section we show how to compute the zeta functions of certain unramified
Galois covers of \dsixcurves.

Let $C$ be a \dsixcurve over a field $K$ of characteristic neither~$2$ nor~$3$,
and suppose all of the Weierstrass points of $C$ are rational over~$K$, so that
$C\cong D_{u,c}$ as in Notation~\ref{N:ucurves} for some $u,c\in K$. We saw
in Section~\ref{sec:family} that there is a $3$-isogeny
$\varphi\colon E_1\to E_2$ of elliptic curves over $K$ such that $\Jac C$ is
isomorphic to $E_1\times E_2$ and such that $C$ has degree-$2$ maps to $E_1$ and
to~$E_2$. Let $P_1$ be the Weierstrass point $(u,0)$ in $C(K)$ and fix double 
covers $\psi_i\colon C\to E_i$ so that $\psi_i(P_1)$ is the identity of~$E_i$.
In this section we will look at the unramified covers of $C$ obtained by pulling
back certain $3$-power isogenies of $E_1$ and $E_2$ via $\psi_1$ and~$\psi_2$.

\begin{notation}
\label{N:pullbacks}
Let $\chi\colon F\to E_i$ be an isogeny whose target is either $E_1$ or $E_2$.
We let $C^{\chi,i}\to C$ denote the pullback of $\chi$ via $\psi_i$. 
\end{notation}

\begin{remark}
\label{R:dualcovers}
Suppose we have analyzed a cover of $C$ obtained from pulling back an isogeny via
the double cover $\psi_1\colon C\to E_1$. If we replace the pair $(u,c)$ with $(3/u, -c)$ 
and follow the same analysis, then Lemma~\ref{lemma:isomorphic D_{u,c}} shows that we will
be looking at the same curve $C$ (up to isomorphism), but now the double cover will
be a map from $C$ to~$E_2$. Thus, it will suffice for us to analyze the covers
of $C$ obtained from $\psi_1$, because the covers coming from $\psi_2$ are
obtained simply by replacing $u$, $r$, and $c$ with $3/u$, $729/r$, and $-c$ 
in all of the formulas we obtain.
\end{remark}

\subsection{The intrinsic triple covers}
\label{ssec:specialtriple}
Among all of the $3$-power isogenies to the curve~$E_1$, there is one that we
have already discussed: namely, the $3$-isogeny
$\phihat\colon E_2\to E_1$. We begin by analyzing the triple cover $C^{\phihat,1}\to C$
obtained by pulling back $\phihat$ via the double cover~$\psi_1$. 
We call this cover, as well as the analogous cover
$C^{\varphi,2}\to C$, the \emph{intrinsic} triple covers of~$C$, because they
arise from data already associated with the curve~$C$.

\begin{lemma}
\label{L:specialPrym}
The Prym variety of the intrinsic triple cover $C^{\phihat,1}\to C$ is isogenous to the
product of the elliptic curve $F^{\phihat,1}$ given by
\[y^2 = x^3 + 81 x^2 + 72 r x + 16 r^2\]
with its twist by $(81-3r)c$.
\end{lemma}

\begin{proof}
Recall that~\eqref{EQ:basicec} gives an equation for $E_1$:
\[
c y^2 =  x^3 + (r-18)x^2 + (81-2r)x + r.
\]
If we take
the equation~\eqref{EQ:ecquotient} for $E_2$ and replace $x$ with $w+r$ and $y$
with~$z$, we get the equation 
\begin{equation}
\label{EQ:E2bis}
cz^2 = w^3 + (r + 81) w^2 + 144 r w + 64 r^2.
\end{equation}
Up to isomorphism, the $3$-isogeny $\phihat$ is then given by
\begin{align*}
x &= \frac{w^3 + 81 w^2 + 72 r w + 16 r^2}{(3w + 4r)^2}\\[1ex]
y &= z \cdot \frac{3 (r - 27) w^2 + c z^2}{(3w + 4r)^3}.
\end{align*}
The cover $\psi_1\colon C\to E_1$ is given by $(x,y)\mapsto (x^2,y)$. The pullback of $\phihat$
along $\psi_1$ is therefore given by~\eqref{EQ:basiccurve} and~\eqref{EQ:E2bis}
together with
\begin{align}
\nonumber x^2 &= \frac{w^3 + 81 w^2 + 72 r w + 16 r^2}{(3w + 4r)^2}\\[1ex]
\label{EQ:ymap} y &= z \cdot \frac{3 (r - 27) w^2 + c z^2}{(3w + 4r)^3}.
\end{align}
If we set $s = (3w + 4r) x$, we find that
\begin{equation}
\label{EQ:special3a} s^2 = w^3 + 81 w^2 + 72 r w + 16 r^2.
\end{equation}
(Note that~\eqref{EQ:special3a} is the curve $F$ from the statement of the lemma.)
We check that the genus-$4$ curve $C^{\phihat,1}$ can be defined by the two equations
\eqref{EQ:E2bis} and \eqref{EQ:special3a} in the variables $w$, $s$, and~$z$,
and the degree-$3$ map from $C^{\phihat,1}$ to $C$ is given by $x = s/(3w+4r)$ and~\eqref{EQ:ymap}.

This model for the curve $C^{\phihat,1}$ shows that it is $V_4$ extension of~$\PP$. 
By~\cite[Theorem~C]{KaniRosen1989}, the Jacobian of $C^{\phihat,1}$ is isogenous to the
product of the elliptic curve~$F^{\phihat,1}$, the elliptic curve given 
by~\eqref{EQ:E2bis} (which is~$E_2$), and the Jacobian of the genus-$2$ 
curve $D$ given by
\begin{equation}
\label{EQ:special3genus2}
cy^2 = (w^3 + (r + 81) w^2 + 144 r w + 64 r^2)(w^3 + 81 w^2 + 72 r w + 16 r^2).
\end{equation}

Define elements $\abar$, $\bbar$, $\cbar$, and $\dbar$ of $K$ by
\[
\abar \colonequals \frac{27}{4(r-27)}, \quad
\bbar \colonequals \frac{-27}{8(r-27)}, \quad
\cbar \colonequals \frac{r-27}{108}, \text{\quad and\ } 
\dbar \colonequals \frac{-(r-27)}{216},
\]
and (reusing the variable $x$) set $x = -3w/(3w + 4r)$, so that 
$w = -4rx/(3x+3)$. Substituting this expression for $w$ 
into~\eqref{EQ:special3genus2} and clearing denominators, we find that $D$ can
be written
\[
c(1-r/27) z^2 = (x^3 + 3\abar x + 2\bbar )(2\dbar x^3 + 3\cbar x^2 + 1)
\]
where
\[z \colonequals y \cdot \frac{(x+1)^3}{2^6 r^2(1-r/27)}.\]
We check that $12 \abar \cbar + 16 \bbar \dbar = 1,$ so by~\cite[\S A.1]{BrokerHoweEtAl2015} we
see that the Jacobian of $D$ is $(3,3)$-isogenous to the product of two
elliptic curves $F_1$ and~$F_2$ given by explicit polynomial equations. We leave
the reader to check that one of these curves is isomorphic to~$E_2$, and the
other can be written
\[
(81 - 3r)c y^2 = x^3 + 81 x^2 + 72 r x + 16 r^2.
\]
This is simply a quadratic twist of the elliptic curve $F^{\phihat,1}$ in the statement of
the lemma. Since the Jacobian of $C$ is isogenous to $E_2^2$, we find that the 
Prym of $C^{\phihat,1}\to C$ is isogenous to the product of the elliptic curve 
$F^{\phihat,1}$ with its twist by $(81-3r)c$.
\end{proof}

\begin{corollary}
\label{C:specialPrym}
The Prym variety of the intrinsic triple cover $C^{\varphi,2}\to C$ is isogenous to the
product of the elliptic curve $F^{\varphi,2}$ given by
\[y^2 = x^3 + r^2 x^2 + 8 r^3 x + 16 r^4\]
with its twist by $r(27-r)c$.
\end{corollary}

\begin{proof}
This follows from Lemma~\ref{L:specialPrym} by applying 
Remark~\ref{R:dualcovers} and rescaling the variables $x$ and $y$.
\end{proof}

\begin{remark}
\label{R:twospecial}
For a generic \dsixcurve $C$, the elliptic curves from Lemma~\ref{L:specialPrym}
are not geometrically isogenous to the elliptic curves from 
Corollary~\ref{C:specialPrym}. To see this, we simply note that for $r=1\in \FF_5$,
the curves in Lemma~\ref{L:specialPrym} are supersingular but the curves
in Corollary~\ref{C:specialPrym} are not.
\end{remark}

\begin{remark}
Another way of characterizing the intrinsic triple covers $C'\to C$ among all
cyclic Galois triple covers is that they are the covers for which the automorphisms $\alpha$ and $\beta$ 
lift to~$C'$. These lifted automorphisms give rise to automorphisms of 
the genus-$2$ curve given by \eqref{EQ:special3genus2}, so this curve is also
a \dsixcurve.

There are two further cyclic Galois triple covers of $C'\to C$ with the property that
the automorphism $\alpha$ lifts to~$C'$. (The automorphism $\beta$ 
lifts to neither; instead, $\beta$ interchanges these two covers.) The Pryms
of these triple covers are isogenous to the square of an elliptic curve with
$j$-invariant $0$, and are therefore not useful for distinguishing curves by
their $L$-functions.
\end{remark}

\subsection{The general elliptic triple covers}
\label{subsec:general triple}
The $3$-isogenies with target curve $E_1$ defined over $\Kbar$ correspond to 
the roots of the $3$-division polynomial $g_3$ for~$E_1$, which is
\begin{align*}
g_3 &= 3x^4 + (4r - 72)x^3 + (486 -12r)x^2 + 12rx + 252r - 6561\\
    &= (x + 3) \bigl[ 3(x-9)^3 + 4 r  (x^2 - 6x + 21)\bigr].
\end{align*} 
Namely, a root $v$ of the division polynomial defines an order-$3$ subgroup
of $E_1(\Kbar)$ and hence a $3$-isogeny from $E_1$ to some elliptic curve~$F$
over~$\Kbar$. The dual of this $3$-isogeny gives an isogeny $F\to E_1$ over $\Kbar$
which we will denote~$\chi_v$.

The root of $g_3$ at $-3$ gives the isogeny $\phihat\colon E_2\to E_1$, so
$C^{\chi_{-3},1}\to C$ is the intrinsic triple cover we analyzed in the
preceding section. In this section we instead consider roots $v$ of the cubic
factor of the division polynomial, so that 
\[
r = \frac{-3(v-9)^3}{4(v^2 - 6v + 21)}.
\]

If $v$ lies in $K$, then $C^{\chi_{v},1}\to C$ is $K$-rational. If 
$L\colonequals K(v)$ is a proper extension of $K$, then $C^{\chi_{v},1}\to C$ 
is rational over~$L$, and if we take the composition of the Galois 
conjugates of the triple cover, we obtain a $3$-power cover of $C$ that is
rational over~$K$. If $L/K$ has degree~$2$, we obtain a degree-$9$ cover of~$C$,
and if $L/K$ has degree~$3$, we obtain a degree-$27$ cover of $C$. Let us
describe how to calculate the Prym varieties of these covers.

\begin{lemma}
\label{L:generalPrym}
Let $v\in \Kbar$ be a root of the division polynomial $g_3$ and let $h$ be the
polynomial 
\[
h\colonequals x^3 + \frac{27(v^2 - 6v + 21)}{4(v - 9)(v + 3)} \, (x + 1)^2.
\]
The Prym variety of the triple cover $C^{\chi_v,1}\to C$ over $\Kbar$
is isogenous to the product of the two elliptic curves
\begin{align}
\label{EQ:gen3b}  3(v-9)cy^2 &= h \quad\text{and}\\
\label{EQ:gen3a} -3(v+3)y^2  &= h.
\end{align}
\end{lemma}

\begin{proof}
Let $E'$ be the curve
\begin{multline}
\label{EQ:general3b}
c z^2 = w^3 - \frac{9}{4}\frac{(3v^3 - 21v^2 + 113v - 159)}{(v^2 - 6v + 21)} w^2 \\
          + 432 \frac{(3v^2 - 10v - 9)}{(v^2 - 6v + 21)} w
          - 20736 \frac{(3v^3 - 17v^2 + 9v + 69)}{(v^2 - 6v + 21)^2}.
\end{multline}
We let the reader check that $\chi_v$ is the isogeny from $E'$ to $E_1$ given by
\begin{align*}
x&= \frac{
w^3 + 9w^2 - \frac{10368}{(v^2 - 6v + 21)} w - \frac{82944(v^2 - 6v + 33)}{(v^2 - 6v + 21)^2}}
{9 \Big( w - \frac{96(v-3)}{(v^2 - 6v + 21)} \Big)^2 }\\
y &= z \cdot \frac{\frac{27(v^3 - 7v^2 - 5v + 75)}{4(v^2 - 6v + 21)}  \Bigl(w - \frac{96}{(v-5)}\Bigr)^2 + cz^2}
        {27 \Bigl(w - \frac{96(v - 3)}{(v^2 - 6v + 21)}\Bigr)^3}\,.
\end{align*}

Analogously to our computation of the Prym of the intrinsic triple covers,
we find that the curve $C^{\chi_{v},1}$ can be defined by the two
equations~\eqref{EQ:general3b} and 
\begin{equation}
\label{EQ:general3a} 
s^2 = w^3 + 9w^2 - \frac{10368}{(v^2 - 6v + 21)} w 
       - \frac{82944(v^2 - 6v + 33)}{(v^2 - 6v + 21)^2},
\end{equation}
so that $C^{\chi_{v},1}$ is a $V_4$ cover of $\PP$, with one intermediate curve
in the extension being $E'$ and one being the elliptic curve $F$ defined 
by~\eqref{EQ:general3a}. The third intermediate curve $D$ is given by setting
$y = sz$ and noting that then $cy^2$ is equal to the product of the right hand
sides of~\eqref{EQ:general3b} and~\eqref{EQ:general3a}. We reuse the 
variable~$x$, and set 
\[
x \colonequals \frac{-(v-5)w/3 + 32}{w - 96(v-3)/(v^2-6v + 21)},
\]
so that
\[ 
w = 32 \cdot \frac{3(v-3)x/(v^2-6v + 21) + 1}{x + (v-5)/3}.
\]
Substituting this expression into the equation for $D$ and clearing
denominators, we find that
\begin{multline*}
c(v-9)^2(v^2-6v+21)^4 (3x + v-5)^6 y^2 \\ 
      = -2^{26} 3^8 r (v+3)^6 (x^3 + 3\abar x + 2\bbar)(2\dbar x^3 + 3\cbar x^2 + 1),
\end{multline*}
where
\begin{align*}
\abar &\colonequals \frac{-v(v^2 - 6v + 21)}{9(v + 3)}\\
\bbar &\colonequals \frac{-(2v^2 - 9v + 27)(v^2 - 6v + 21)}{54(v + 3)}\\
\cbar &\colonequals \frac{-9(v - 5)}{4(v - 9)(v^2 - 6v + 21)}\\
\dbar &\colonequals \frac{27}{8(v - 9)(v^2 - 6v + 21)}\, .
\end{align*}
By scaling $y$, we can rewrite the equation for $D$ as
\[
-cr y^2  = (x^3 + 3\abar x + 2\bbar)(2\dbar x^3 + 3\cbar x^2 + 1).
\]
We check that $12 \abar \cbar + 16 \bbar \dbar = 1,$ so once again 
by~\cite[\S A.1]{BrokerHoweEtAl2015} we see that the Jacobian of $D$ is
$(3,3)$-isogenous to the product of two elliptic curves $F_1$ and~$F_2$ given by
explicit polynomial equations $y^2 = f_1$ and $y^2=f_2$, with the $f_i$ as
in~\cite[\S A.1]{BrokerHoweEtAl2015}. In particular, $F_1$ can be written
\[
-cr y^2 = \Bigl(x - \frac{(v-9)^2}{12(v^2 + 6v + 9)}\Bigr)^3 
       + \frac{3(v+3)}{4(v-9)}\Bigl(x + \frac{(v-9)^2}{36(v^2+6v+9)}\Bigr)^2,
\]
and we check that this curve is $3$-isogenous to~$E_1$.

We check that the curve $F_2$ is isomorphic to the curve defined by~\eqref{EQ:gen3b},
while the curve given by~\eqref{EQ:general3a} is isomorphic to $3r(v+3)y^2 = h$.
However, since our curve $C$ has a rational Weierstrass point, we see from
Lemma~\ref{L:rationalW} that $-r$ is a square, so in fact the
curve~\eqref{EQ:general3a} is isomorphic to the curve defined by~\eqref{EQ:gen3a}.

All told, we see that the Jacobian of $C^{\chi_{v},1}$ is isogenous to the 
product of $E_1^2$ and the two elliptic curves~\eqref{EQ:gen3b} 
and~\eqref{EQ:gen3a}. Since the Jacobian of $C$ is isogenous to $E_1^2$, the 
Prym variety of the triple cover $C^{\chi_{v},1}\to C$ is isogenous to the 
product of~\eqref{EQ:gen3b} and~\eqref{EQ:gen3a}.
\end{proof}

The cover $C^{\chi_{v},1}\to C$ is only defined over the extension $L/K$, and 
what we would really like to know is the isogeny class of the Prym of the extension 
$\Ctilde\to C$ obtained as the composition of $C^{\chi_{v},1}\to C$ with all of its Galois
conjugates.

The maps from $C^{\chi_{v},1}$ to~\eqref{EQ:gen3b} and~\eqref{EQ:gen3a} are defined over~$L$,
and the Galois conjugates of these maps give maps from the conjugates of $C^{\chi_{v},1}$ to
the conjugates of the elliptic curves. Therefore, all of these maps fit together
to produce an isogeny from the Prym of $\Ctilde\to C$ to the product of the 
restrictions of scalars of~\eqref{EQ:gen3b} and of~\eqref{EQ:gen3a} from $L$
to~$K$. 

\begin{remark}
\label{R:tripleWeil}
If we are working over a finite field, the preceding observation gives us an 
efficient way of computing the Weil polynomial of the Prym variety of the 
cover $\Ctilde\to C$. If the degree of $L/K$ is~$e$, and if the Weil polynomial
of~\eqref{EQ:gen3b} is $x^2 - tx + q^e$, then the Weil polynomial of the
restriction of scalars of~\eqref{EQ:gen3b} is $x^{2e} - t x^e + q^e$, and the
analogous statement is true for the restriction of scalars of~\eqref{EQ:gen3a}.
\end{remark}

\section{Initial heuristics and data} 
\label{sec:2coverheuristics}

In this section we give some initial heuristics for how often we expect to find
pairs of \dsixcurves that are doubly isogenous to one another. To begin with,
we observe that we have good estimates for how many pairs of \dsixcurves have
isogenous Jacobians. 

\begin{proposition}
For $q$ be a prime power coprime to $6$, let $P(q)$ denote the number of pairs of \dsixcurves over $\FF_q$ whose
Jacobians are isogenous to one another. Then there are positive constants $d_1$ 
and $d_2$ such that for all sufficiently large prime powers~$q$ coprime to $6$, we have
\[
d_1q^{3/2}  \le P(q) \le d_2q^{3/2}(\log q)^2(\log{\log q})^4.
\]
\end{proposition}

\begin{proof}
The analysis of $D_4$ curves from \cite[\S 4]{ArulEtAl2024} carries through 
\emph{mutatis mutandis} to our current situation. In particular, the argument
that proves \cite[Theorem 4.6]{ArulEtAl2024} proves this proposition as well. 
As the Jacobian of a \dsixcurve is isogenous to the square of an elliptic curve,
this reduces to comparing the number of elliptic curves over $\FF_q$ with the 
number of isogeny classes.  The lower bound is an easy argument, while the 
upper bound uses non-trivial results on the maximum size of isogeny classes 
of elliptic curves over~$\FF_q$.
\end{proof}

Remember from Definition~\ref{def:doubly isogenous} that two \dsixcurves $C_1$ 
and $C_2$ are doubly isogenous  if $\Jac C_1$ is isogenous to $\Jac C_2$ and 
$\Jac C_1^{[2]}$ is isogenous to $\Jac C_2^{[2]}$.  For this to happen, there 
must be multiple simultaneous isogenies among the elliptic curves that appear as 
factors of the Jacobians of $C_1^{[2]}$ and $C_2^{[2]}$. To simplify matters, 
for our heuristics we will consider \emph{geometric} isogenies among these 
curves instead of $\FF_q$-rational isogenies, with the expectation that up to
constant factors the result of the analysis would be the same. That is, if we
view the elliptic curves that arise as factors of our Jacobians to be somehow
drawn at random from the elliptic curves over $\FF_q$, then the probability that 
two geometrically isogenous curves are actually isogenous over $\FF_q$ is
approximately~$1/2$.

\begin{example}
\label{example:doubleisogenypatterns}
Let $D_{u,c}$ be a curve in the \dsixfamily given in Notation~\ref{N:ucurves},
and consider the maximal unramified elementary abelian $2$-cover 
${\pi\colon D^{[2]}_{u,c}\to D_{u,c}}$ described in 
Section~\ref{sec:covers2}. The Jacobian of $D_{u,c}$ is isogenous to the 
square of the elliptic curve $E_{r,c}$ in 
Corollary~\ref{P:baseE}, where $r = -u^2(u^2-9)^2/(u^2-1)^2$, 
and the Prym variety of the cover $\pi$ is also
isogenous to a product of elliptic curves.  As noted in Proposition~\ref{P:2decomp},
up to geometric isomorphism there
are four elliptic curves appearing in this Prym, one with multiplicity $6$ and 
three with multiplicity $3$; see Table~\ref{table:lambdas}. We take $E^{(6)}$, 
$E^{(\threea)}$, $E^{(\threeb)}$, and $E^{(\threec)}$ to be the curves with the $\lambda$-invariants 
given in the table, setting $e=1$ for each. These curves depend on the 
parameter~$u$ associated to the \dsixcurve we started with.

Now consider what 
it means for two members of the family $D_{u_1,c_1}, D_{u_2,c_2}$ to be (geometrically)
doubly isogenous. First, $E_{r_1,c_1}$ and $E_{r_2,c_2}$ must be isogenous. 
Furthermore, there must be isogenies
between $E^{(6)}, E^{(\threea)}, E^{(\threeb)}, E^{(\threec)}$ for the first curve 
and for the second. One possibility is to have isogenies
$E^{(6)}_{u_1} \sim E^{(6)}_{u_2}$, $E^{(\threea)}_{u_1} \sim E^{(\threea)}_{u_2}$, and
so on.  Other possibilities interchange curves, for example having isogenies 
$E^{(6)}_{r_1} \sim E^{(\threea)}_{r_2} \sim E^{(\threeb)}_{r_2}$, 
$E^{(\threea)}_{r_1} \sim E^{(\threeb)}_{r_1} \sim E^{(6)}_{r_2}$, and 
$E^{(\threec)}_{r_1} \sim E^{(\threec)}_{r_2}$.  There are a variety of
configurations, all of which require (at least) four isogenies between the 
elliptic curves in the Pryms in addition to the isogeny between the base 
elliptic curves.
\end{example}

When looking for isogenies, in the absence of additional information our default
heuristic is the following (cf. \cite[Lemma~5.3]{ArulEtAl2024} but note that the
situation in that paper leads to six elliptic curves). 

\begin{heuristic} \label{H:default}
The elliptic curve $E_{r,c}$ and the four elliptic curves $E^{(6)}_u$, 
$E^{(\threea)}_u$, $E^{(\threeb)}_u$, $E^{(\threec)}_u$ that appear in the Jacobian
of $D_{u,c}^{[2]}$ behave as independent random elliptic curves over $\FF_q$.
\end{heuristic}

In particular, this predicts that the probability that two of these elliptic 
curves are isogenous for different values of the parameter $u$ is $\asymp q^{-1/2}$.

\begin{corollary}
Assuming Heuristic~\textup{\ref{H:default}}, the expected number of pairs of
doubly isogenous \dsixcurves over $\FF_q$ is $\asymp q^{-1/2}$.  In particular,
we expect to see infinitely many such doubly isogenous pairs over the totality
of all finite fields.
\end{corollary}

\begin{proof}
There are $\asymp q^2$ pairs of \dsixcurves and to be doubly isogenous means 
that $5$ pairs of elliptic curves are isogenous, so the expected
number of doubly isogenous pairs is $\asymp q^2 (q^{-1/2})^5 = q^{-1/2}$.  The
last statement follows as $\sum_q q^{-1/2}$ is divergent.
\end{proof}

Now consider the probability of encountering a doubly isogenous
pair for which the Jacobians of the intrinsic triple covers from 
Section~\ref{ssec:specialtriple} are also isogenous. There is an obvious
generalization of Heuristic~\ref{H:default} to this situation:

\begin{heuristic} \label{H:default-special}
The five elliptic curves in Heuristic~\textup{\ref{H:default}}, the
elliptic curve $F^{\phihat,1}$ from Lemma~\textup{\ref{L:specialPrym}}, and
the elliptic curve $F^{\varphi,2}$ from Corollary~\textup{\ref{C:specialPrym}}
behave as independent random elliptic curves over~$\FF_q$.
\end{heuristic}

This heuristic suggests that, absent additional phenomena, the expected number
of doubly isogenous pairs over $\FF_q$ for which the Pryms of the intrinsic 
triple covers are also isogenous will be ${\asymp q^{-3/2}}$. Since
$\sum_q q^{-3/2}$ is a convergent series, we expect to see only finitely many
such pairs over all finite fields.
 
We will also want to consider the likelihood of finding doubly isogenous pairs
for which the Pryms of higher degree unramified covers are isogenous. For
example, given a \dsixcurve $C$ with Jacobian isomorphic to $E_1\times E_2$ as 
in Proposition~\ref{T:S3surfaces}, we might consider the cover $D\to C$ obtained by
pulling back the multiplication-by-$3$ map on $E_1$ via the double cover 
$C\to E_1$. Geometrically, this cover is abelian with group $(\ZZ/3\ZZ)^2$, and
its Prym decomposes up to isogeny as the product of the Pryms of the geometric 
triple covers $C^{\chi_v,1}\to C$ discussed in Section~\ref{sec:triplecovers}. 
Over~$\FF_q$, the Prym of this cover decomposes as the product of the 
restrictions of scalars discussed in Remark~\ref{R:tripleWeil}.

Suppose, for example, that we have two curves $C_1$ and $C_2$ for which the
Prym of this cover is isogenous to the square of the elliptic curve coming from
the intrinsic triple cover of Section~\ref{ssec:specialtriple} and the square of
the restriction of scalars of an elliptic curve defined over $\FF_{q^3}$ of the
type considered in Section~\ref{subsec:general triple}. What is the probability 
that the three-dimensional factors are isogenous? The spirit of our heuristic is
that the probability should be that of having two random elliptic curves over 
$\FF_{q^3}$ being isogenous, which would be $\asymp q^{-3/2}$. The fact we
observe here is that this probability is the same order of magnitude as that of
having three pairs of elliptic curves over $\FF_q$ be isogenous.

The point is that we can use our estimate that two Pryms will be isogenous with
probability equal to $q^{-m/2}$, where $m$ is the number of distinct elliptic
factors of the Prym, \emph{even when these factors are only defined 
geometrically}.

So, if instead of only considering the intrinsic triple covers of 
Section~\ref{ssec:specialtriple}, we look at the covers obtained from the
multiplication-by-$3$ maps on the two elliptic factors of the Jacobians of 
curves in our \dsixfamily, we find that there are a total of eight more
elliptic curves (in addition to the five from Corollary~\ref{P:baseE} and 
Table~\ref{table:lambdas}). 
An extension of Heuristic~\ref{H:default} is:

\begin{heuristic}
The five elliptic curves in Heuristic~\textup{\ref{H:default}}
and the additional eight elliptic curves appearing in the Pryms of the triple
covers $C^{\chi_v,i}\to C$ from Section~\textup{\ref{sec:triplecovers}} behave as 
independent random elliptic curve over~$\FF_q$.
\end{heuristic}

This predicts that the probability
of two \dsixcurves being doubly isogenous and having these $3$-power covers 
also being isogenous is $\asymp q^{-13/2}$. The expected number of such pairs
over $\FF_q$ is therefore $\asymp q^{-9/2}$.  We expect finitely many such pairs.

Finally, we can consider the $26$ orbits of elliptic curves obtained
from unramified fourfold covers of~$C$.

\begin{heuristic}
The five elliptic curves in Heuristic~\textup{\ref{H:default}}
and the additional $26$ curves from Lemma~\textup{\ref{L:4action}}
behave as independent random elliptic curve over $\FF_q$.
\end{heuristic}

If we ask for a pair of doubly isogenous curves 
with the $4$-power covers of Section~\ref{sec:covers4} being isogenous, this heuristic predicts this happens with probability of $\asymp q^{-31/2}$. The
expected number of such pairs over $\FF_q$ is then $\asymp q^{-27/2}$.   We expect finitely many such pairs.

To check our heuristics, we gathered data. For each $n$ from $12$ up to $30$, we
looked at the $1024$ primes $p$ closest to $2^n$, and for each such $p$ we
enumerated all pairs $\{D_{u_1,c_1}, D_{u_2,c_2}\}$ of doubly isogenous curves 
in the family $D_{u,c}$. We separated the pairs into two groups: those where the
curves are twists of one another, and those where the curves are geometrically
non-isomorphic. (This is natural to do, because our heuristic about the 
independence of the elliptic curves associated to two \dsixcurves clearly fails
when the \dsixcurves are twists of one another.) For each~$n$, we calculated 
the total number of doubly isogenous pairs (twists and nontwists) for these 
$1024$ primes~$p$. The data is presented in the second and sixth columns of Table~\ref{table:doublyplus}.

The code we used to compute these values is available in the GitHub 
repository for this paper~\cite{DoublyIsogenousRepo}. We ran the code on several systems:
\begin{itemize}
\item \babbage, a 64-core Intel Xeon E7-8867v3 CPU running at 3.3GHz with 1TB 
      memory;
\item \hensel, a 64-core AMD Ryzen Threadripper PRO 3995WX CPU running at 2.7GHz
      with 1TB memory; and
\item \winifred, a 10-core Apple M1 Max running at 3.2 GHz with 64 GB memory.
\end{itemize}
The computations took a total of 539 CPU-days on \hensel, 378 CPU-days on 
\babbage, and 11 CPU-days on \winifred.

\begin{table}[t]
\centering
\caption{Data for doubly isogenous \dsixcurves for which the Pryms of further
covers are also isogenous. For each $n$, the second column lists the total
number of (unordered) pairs of doubly isogenous \dsixcurves over $\Fq$ that are
twists of one another, for the $1024$ primes $q$ closest to $2^n$. The third
column lists the number of such curves that also have isogenous Pryms from the
intrinsic triple covers of Section~\ref{sec:triplecovers}, the
fourth column lists the number with isogenous Pryms from the pullbacks of the
multiplication-by-$3$ maps of the underlying elliptic curves, and the fifth
column lists the number with isogenous Pryms from the pullbacks of the
multiplication-by-$4$ map on the Jacobian. The sixth through ninth columns give
the same information for pairs of curves that are not twists of one another.}
\label{table:doublyplus}
\begin{tabular}{ccccccccccc}
\toprule
\hbox to 2em{}&&\multicolumn{4}{c}{Twists}&\hbox to 0.5em{}&\multicolumn{4}{c}{Nontwists}\\
\cmidrule(lr){3-6}\cmidrule(lr){8-11}
$n$ && $2$ & $3$a & $3$b & $4$ && $2$ & $3$a & $3$b & $4$\\
\midrule
$12$ && $267$ &    $25$ & $0$ & $189$ &&     $2568$ & $150$ & $2$ & $0$ \\
$13$ && $285$ &    $24$ & $0$ & $228$ &&     $2294$ & $145$ & $3$ & $0$ \\
$14$ && $295$ &    $18$ & $1$ & $256$ &&     $1749$ & $135$ & $0$ & $0$ \\
$15$ && $283$ &    $12$ & $0$ & $249$ &&     $1428$ & $135$ & $0$ & $0$ \\
$16$ && $250$ & \pz $6$ & $0$ & $233$ &&     $1064$ & $133$ & $0$ & $0$ \\
$17$ && $260$ & \pz $5$ & $0$ & $242$ && \pz  $860$ & $133$ & $0$ & $0$ \\
$18$ && $269$ & \pz $6$ & $0$ & $257$ && \pz  $631$ & $117$ & $0$ & $0$ \\
$19$ && $271$ & \pz $3$ & $0$ & $261$ && \pz  $571$ & $151$ & $0$ & $0$ \\
$20$ && $262$ & \pz $1$ & $0$ & $255$ && \pz  $422$ & $135$ & $0$ & $0$ \\
$21$ && $243$ & \pz $2$ & $0$ & $242$ && \pz  $328$ & $106$ & $0$ & $0$ \\
$22$ && $264$ & \pz $0$ & $0$ & $261$ && \pz  $282$ & $114$ & $0$ & $0$ \\
$23$ && $268$ & \pz $0$ & $0$ & $265$ && \pz  $270$ & $129$ & $0$ & $0$ \\
$24$ && $291$ & \pz $0$ & $0$ & $291$ && \pz  $236$ & $131$ & $0$ & $0$ \\
$25$ && $236$ & \pz $0$ & $0$ & $236$ && \pz  $222$ & $132$ & $0$ & $0$ \\
$26$ && $213$ & \pz $0$ & $0$ & $212$ && \pz  $232$ & $144$ & $0$ & $0$ \\
$27$ && $260$ & \pz $0$ & $0$ & $260$ && \pz  $194$ & $132$ & $0$ & $0$ \\
$28$ && $252$ & \pz $0$ & $0$ & $251$ && \pz  $210$ & $153$ & $0$ & $0$ \\
$29$ && $286$ & \pz $0$ & $0$ & $286$ && \pz  $150$ & $107$ & $0$ & $0$ \\
$30$ && $245$ & \pz $0$ & $0$ & $245$ && \pz  $183$ & $124$ & $0$ & $0$ \\
\bottomrule
\end{tabular}
\end{table}

The data seems to suggest that as $n$ increases the number of doubly isogenous twists is
remaining roughly constant, and that the number of doubly isogenous nontwists
is perhaps tending towards a nonzero value. In the following section we explain
why this is so, and we suggest revisions to our heuristics.

In any case, it appears that doubly isogenous pairs in our family of \dsixcurves
are not uncommon. Therefore, if we want to distinguish curves from one another by isogeny
criteria, we will need to use more than just the pullback of the 
multiplication-by-$2$ map on the Jacobian. To see whether isogeny classes of 
further covers actually help to distinguish curves, we look at the same set of
primes as before, and we check to see how many pairs of doubly isogenous curves
also have isogenous Pryms coming from the two intrinsic triple covers of the
underlying elliptic factors, how many pairs have isogenous Pryms coming from the
pullbacks of multiplication-by-$3$ on the underlying elliptic factors, and how
many pairs have isogenous Pryms coming from the pullback of the 
multiplication-by-$4$ map on the Jacobian. The data is presented in the remaining columns of 
Table~\ref{table:doublyplus}.

We see that for doubly isogenous pairs of \dsixcurves that are twists of one
another, the intrinsic triple covers tend to be able to distinguish the curves,
and the triple covers from the elliptic factors distinguish all but one pair of
curves in our data set. However, the Pryms of the covers coming from the
multiplication-by-$4$ map on the Jacobian do not seem to help much in
distinguishing the curves from one another.

On the other hand, for doubly isogenous pairs of \dsixcurves that are
\emph{not} twists of one another, the intrinsic triple covers do not seem to be
sufficient to distinguish curves. The general triple covers coming from the
elliptic factors do a much better job at distinguishing these curves from one 
another, and the Pryms coming from the multiplication-by-$4$ maps distinguish
all of the pairs of such curves in our data set.

\begin{remark} \label{R:noisogenies}
The code used to gather this data can also easily show that there exist
$u,v \in \QQ$ such that the ten elliptic curves appearing in the decomposition
of $\Jac D_{u,1}$ and $\Jac D_{v,1}$ are all pairwise non-isogenous
over~$\Qbar$.  For example, if we take $u=34$ and $v=57$, we find that the
reductions of the ten elliptic curves to $\FF_{509}$ are all ordinary, so the
geometric endomorphism rings of the reductions are isomorphic to their 
endomorphism rings over~$\FF_{509}$. We compute that the discriminants of these
ten endomorphism rings are all distinct in $\QQ^\times/\QQ^{\times 2}$. But the
class in $\QQ^\times/\QQ^{\times 2}$ of the discriminant of the endomorphism 
ring of an ordinary elliptic curve over the algebraic closure of a finite field
is an isogeny invariant. This shows that the reductions modulo $509$ of the ten
curves are pairwise geometrically non-isogenous, so the same is true of the
curves themselves.
\end{remark}

\section{Unexpected sources for doubly isogenous curves}
\label{sec:extraordinary}

In this section we present two constructions that explain the overabundance of 
doubly isogenous curves we observed in the preceding section.

\subsection{Doubly isogenous twists}
\label{subsec:twists}

As we can see, the number of doubly isogenous twists in Table~\ref{table:doublyplus}
is roughly constant for all~$n$. This can be explained as follows. 

\begin{theorem}
\label{T:easytwists}
Let $p$ be a prime with $p\equiv 3 \bmod 4$ and let $u$ be an element of $\FF_p$ 
that is not an element of $U_\bad$ \textup(see 
Notation~\textup{\ref{N:ucurves}}\textup). Suppose that
\[
\frac{(u-1)^3(u+3)}{(u+1)^3(u-3)} 
    \text{\quad and\quad} 
\frac{4u}{(u-1)(u+3)}
\]
are supersingular $\lambda$-invariants, and that $u(u-1)(u+3)$ and 
$-u(u+1)(u-3)$ are nonsquares in $\FF_p^\times$. Then the curve $D_{u,1}$ and its quadratic twist
are doubly isogenous. Furthermore, the Pryms of the pullbacks of the 
multiplication-by-$4$ maps on the Jacobians of $D_{u,1}$ and its twist are isogenous.
\end{theorem}

\begin{proof}
Let $c$ be a nonsquare in $\FF_p$, so that our goal is to show that $D_{u,1}$
and $D_{u,c}$ are doubly isogenous.

Let $r=-u^2(u^2-9)^2/(u^2-1)^2$ and let
\[
\lambda=\frac{(u-1)^3(u+3)}{(u+1)^3(u-3)} .
\]
We check that the $j$-invariant of the curve $y^2 = x(x-1)(x-\lambda)$ is equal
to that of the elliptic curve $E_{r,1}$ from Corollary~\ref{P:baseE}. Thus, the 
first hypothesis of the theorem implies that $E_{r,1}$ is supersingular, and 
hence has trace~$0$ because we are working over~$\FF_p$. Therefore $E_{r,1}$ is 
isogenous to its twist $E_{r,c}$, and it follows that the Jacobian of $D_{u,1}$
is isogenous to the Jacobian of $D_{u,c}$.

Next we will show that for each of the orbits in Table~\ref{table:lambdas}, the 
multiset of traces of the elliptic curves in that orbit for $D_{u,1}$ is equal
to the corresponding multiset for $D_{u,c}$. We start with Orbit~6, for which we
have
\[
\lambda = \frac{4u}{(u-1)(u+3)}.
\]
The second hypothesis of the theorem is that this $\lambda$-invariant is
supersingular, so all of the elliptic curves in Orbit~6 have trace~$0$, both for
$D_{u,1}$ and $D_{u,c}$.

Now consider Orbit~\threeb. There are three values of $e$ in 
Table~\ref{table:lambdas} for the elliptic curves in Orbit~\threea. The first value
does not depend on~$c$, so we get the same trace for $D_{u,1}$ and 
for~$D_{u,c}$. Since $u(u-1)(u+3)$ is a nonsquare by hypothesis, the product of
the second and third values of $e$ is a nonsquare, so the second and third
values of $e$ lie in different square classes. It follows that the two
associated elliptic curves have opposite traces, and this is true whether we are
looking at $D_{u,1}$ or $D_{u,c}$. Therefore, the multiset of traces of the 
elliptic curves in Orbit~\threeb\ for $D_{u,1}$ is equal to the corresponding multiset
for $D_{u,c}$.

Since $-u(u+1)(u-3)$ is also assumed to be nonsquare, the same argument shows
that the multiset of traces of the elliptic curves in Orbit~\threec\ for $D_{u,1}$ is
equal to the corresponding multiset for $D_{u,c}$. Finally, since $u(u-1)(u+3)$
and $-u(u+1)(u-3)$ and $-1$ are all nonsquares, so is $(u-1)(u+1)(u-3)(u+3)$, 
which is their product divided by a square. Thus, the same argument shows that
the multiset of traces of the elliptic curves in Orbit~\threea\ for $D_{u,1}$ is equal
to the corresponding multiset for~$D_{u,c}$. It follows that $D_{u,1}$ and 
$D_{u,c}$ are doubly isogenous.

We noted in Section~\ref{sec:covers4} that the Weil polynomial of the Prym of
the cover ${D_{u,1}^{[4]}\to D_{u,1}^{[2]}}$ can be computed as the product of 
Weil polynomials coming from elliptic factors. For each geometric elliptic 
factor, we computed a set of expression in $u$ and $c$ that must be squares in
order for the geometric factor to be rational over~$\FF_q$. If the elliptic
factor is rational over $\FF_q$, its Weil polynomial over $\FF_q$ is a factor of
the Weil polynomial of the Prym; if the elliptic factor is only rational over 
$\FF_{q^2}$, then the Weil polynomial of its restrictions of scalars, from 
$\FF_{q^2}$ to $\FF_q$, is a factor of the Weil polynomial of the Prym. We check
that if $-1$, $u(u-1)(u+3)$, and $-u(u-1)(u+3)$ are all nonsquare, then none of 
the elliptic factors are rational over~$\FF_q$. But over $\FF_{q^2}$, the 
elliptic factors for $D_{u,1}$ and for $D_{u,c}$ are identical, because $c$ is a
square in $\FF_{q^2}$. Therefore, the Pryms of the pullbacks of the 
multiplication-by-$4$ map on the Jacobians of $D_{u,1}$ and its twist are not
just isogenous to one another, they are isomorphic to one another.
\end{proof}

If we view ``being supersingular'' as an event of probability $\asymp p^{-1/2}$
for elliptic curves over $\FF_p$, and if we view the supersingularity of the two
$\lambda$s in the theorem as being events independent of one another and
independent of the ``nonsquare'' conditions in the theorem, then we might view 
the conditions of the theorem being satisfied as an event with probability
$\asymp p^{-1}$. Since there are at least $p/12-1$ possible values of~$u$ (up to
the equivalence of Lemma~\ref{lemma:isomorphic D_{u,c}}), it is reasonable to
think that the expected number of doubly isogenous twists over $\FF_p$ is
positive. In fact, we can prove that this expected value is positive under certain congruence conditions.

\begin{corollary}
\label{C:easytwists}
Let $p$ be a prime with $p\equiv 11\bmod 120$ or $p\equiv 59\bmod 120$, let 
$s,t\in \FF_p$ satisfy $s^2 = 3$ and $t^2 = 5$, and let $u = s(1+t)/2$, so that
$u^4 - 9u^2 + 9 = 0$. Then the curve $D_{u,1}$ and its quadratic twist are 
doubly isogenous, and they cannot be distinguished by the Pryms of the pullbacks
of the multiplication-by-$4$ maps on the Jacobian.
\end{corollary}

\begin{proof}
We note that the congruence conditions on $p$ imply that $-1$ is not a square
modulo~$p$ but $-2$, $3$, and $5$ are squares modulo~$p$, so in particular 
there do exist elements $s$ and $t$ as in the statement of the corollary. We
claim that the conditions of Theorem~\ref{T:easytwists} hold.

Set 
\[
\lambda \colonequals \frac{(u-1)^3(u+3)}{(u+1)^3(u-3)}.
\]
We compute that the corresponding $j$-invariant is a root of 
\[
x^2 - 37018076625x + 153173312762625,
\]
which is the Hilbert class polynomial for the discriminant $-60$. This means
that the elliptic curve corresponding to $\lambda$ has CM by $\ZZ[\sqrt{-15}]$,
the imaginary quadratic order of discriminant~$-60$. But since
$-60$ is not a square modulo~$p$, a curve over $\FF_p$ with CM by $\ZZ[\sqrt{-15}]$ is 
supersingular.

Similarly, if we set 
\[
\lambda \colonequals \frac{4u}{(u-1)(u-3)},
\]
we find that the corresponding $j$-invariants is $54000$, which is the 
$j$-invariant of an elliptic curve with CM by $\ZZ[\sqrt{-3}]$. Since $-3$ is not a square
modulo~$p$, this elliptic curve is supersingular as well.

We note that
\[
2u(u-1)(u+3) = (s+3)^2(t+1)^2/4,
\]
so that $2u(u-1)(u+3)$ is a square. Since $2$ is a nonsquare, we see that 
$u(u-1)(u+3)$ is also a nonsquare. Similarly, 
\[
-2u(u+1)(u-3) = (s-3)^2(t+1)^2/4,
\]
so $-u(u+1)(u-3)$ is a nonsquare.

Thus, the hypotheses of Theorem~\ref{T:easytwists} hold, so $D_{u,1}$ and its 
quadratic twist are doubly isogenous.
\end{proof}

We can ask how often doubly isogenous twists are explained by 
Theorem~\ref{T:easytwists}. In Table~\ref{table:easytwists} we go through the
data on doubly isogenous twists from Table~\ref{table:doublyplus}, and break down 
the number of twists into the number that are explained by 
Theorem~\ref{T:easytwists} and the number that are not.

\begin{table}[t]
\centering
\caption{Data for doubly isogenous curves that are twists of one another.
For each $n$, the second column lists the total number of (unordered) pairs of 
doubly isogenous curves over $\Fq$ that are twists of one another, for the 
$1024$ primes $q$ closest to $2^n$. The third column lists the number of these
examples that satisfy the hypotheses of Theorem~\ref{T:easytwists}, and the
fourth column lists the number that do not satisfy those hypotheses.}
\label{table:easytwists}
\begin{tabular}{ccrccccrc}
\toprule
& & \multicolumn{2}{c}{Theorem~\ref{T:easytwists}?}&\hbox to 2em{}&& & \multicolumn{2}{c}{Theorem~\ref{T:easytwists}?}\\
\cmidrule(lr){3-4}\cmidrule(lr){8-9}
$n$& Twists & \ \ Yes & No && $n$ & Twists & \ \ Yes & No \\
\cmidrule(lr){1-4}\cmidrule(lr){6-9}
$12$ & $267$ & $188$ &    $79$ && $22$ & $264$ & $261$ & $3$ \\
$13$ & $285$ & $228$ &    $57$ && $23$ & $268$ & $265$ & $3$ \\ 
$14$ & $295$ & $256$ &    $39$ && $24$ & $291$ & $291$ & $0$ \\ 
$15$ & $283$ & $249$ &    $34$ && $25$ & $236$ & $236$ & $0$ \\ 
$16$ & $250$ & $233$ &    $17$ && $26$ & $213$ & $212$ & $1$ \\ 
$17$ & $260$ & $242$ &    $18$ && $27$ & $260$ & $260$ & $0$ \\ 
$18$ & $269$ & $257$ &    $12$ && $28$ & $252$ & $251$ & $1$ \\ 
$19$ & $271$ & $261$ &    $10$ && $29$ & $286$ & $286$ & $0$ \\ 
$20$ & $262$ & $255$ & \pz $7$ && $30$ & $245$ & $245$ & $0$ \\ 
$21$ & $243$ & $242$ & \pz $1$ &&    &     &     &   \\
\bottomrule
\end{tabular}
\end{table}

Of course, Theorem~\ref{T:easytwists} gives just one possible way for a \dsixcurve
over a finite prime field to be doubly isogenous to its twist. We can 
show, however, that every \dsixcurve over $\FF_p$ that is doubly isogenous to 
its twist heuristically requires an event of probability~$\asymp 1/p$. Let us
see why:

First of all, for a \dsixcurve $D_{u,1}$ to be isogenous to its twist requires
that the underlying elliptic curve $E_{r,1}$ be isogenous to \emph{its} twist,
and in particular $E_{r,1}$ must have trace~$0$, and so must be supersingular.
Up to logarithmic factors, the probability that a random elliptic curve over 
$\FF_p$ is supersingular is $1/\sqrt{p}$.

Suppose all of the Orbits~6, \threea, \threeb, and~\threec\ have 
distinct traces (up to sign) and none of them have trace~$0$. Then the set of
traces in each of the orbits for $D_{u,1}$ must match those of its twist, and
by looking at the values of $e$ for Orbit~\threea\ in Table~\ref{table:lambdas} we 
see that this means that $(u-3)(u+1)(u^2+3)$ and 
$(u-1)(u+3)(u^2+3)$ lie in different square classes. That means that 
$d_\threea\colonequals (u-3)(u-1)(u+1)(u+3)$ must be a nonsquare.

Similarly, by looking at Orbit~\threeb, we find that 
$d_\threeb\colonequals u(u-1)(u+3)$ must be a nonsquare, and by looking at 
Orbit~\threec, we find that $d_\threec\colonequals -u(u-3)(u+1)$ must be a
nonsquare. Therefore, the product $d_\threea d_\threeb d_\threec$ is a 
nonsquare. But this product is $-1$ times a square, so we see that $-1$ must be
a nonsquare in~$\FF_p$.

On the other hand, by looking at the values of $e$ for Orbit~6, we see that for
$D_{u,1}$ and its twist to have the same Orbit~6 traces, exactly two of the
following values must be squares:
\begin{align*}
-u(u-3)(u-1)(u+1)(u+3)(u^2+3), && u(u^2+3),          \\
(u-3)(u+1)(u^2+3),             && (u-1)(u+3)(u^2+3).
\end{align*}
Therefore the product of these values must be a square. But the product is $-1$
times a square, so $-1$ must be a square in $\FF_p$. 

This contradiction shows that if $D_{u,1}$ is doubly isogenous to its twist, 
then either one of its orbits consists of supersingular curves, or two of its
orbits have the same trace (up to sign). Up to logarithmic factors, both of
these possibilities have heuristic probability $\asymp 1/\sqrt{p}$. Thus, 
heuristically, the probability that a random \dsixcurve over $\FF_p$ is doubly
isogenous to its twist is $\asymp 1/p$.

\begin{remark}
We can also ask about curves $D_{u,1}$ over non prime finite fields that are
doubly isogenous to their quadratic twists. It turns out that if 
$D_{u,1}/\FF_{p^f}$ is such a curve, then $f$ is odd and $u$ lies in~$\FF_p$.

To see this, we recall from our argument above that if $D_{u,1}$ is such a curve
then the underlying elliptic curve $E_{r,1}$ is isogenous to its twist, and
therefore has trace $0$ and is supersingular. Now, every supersingular elliptic
curve in characteristic $p$ has a model over $\FF_{p^2}$ with all of its 
order-$3$ subgroups rational over $\FF_{p^2}$ and with all of its $2$-torsion
points rational over $\FF_{p^2}$. It follows from Lemma~\ref{L:X0(3)} that if 
$j\in\FF_{p^2}$ is a supersingular $j$-invariant, then all of the values of $r$
in $\Fbar_{p}$ satisfying $-(r-3)^3(r-27)/r = j$ lie in $\FF_{p^2}$, because 
that expression in $r$ is the $j$-invariant of the curve in the lemma. 
Therefore, the values of $r\in\Fbar_p$ such that $C_{r,1}$ is supersingular all
lie in $\FF_{p^2}$. Furthermore, since all of the $2$-torsion of the Jacobian of
$C_{r,1}$ is rational over $\FF_{p^2}$, the values of $u$ corresponding to $r$
are all elements of~$\FF_{p^2}$ as well.

One can show that for all such $r$, the curve~\eqref{EQ:basicec} is isomorphic
to the quadratic twist by $(r-3)(r^2-18r-27)$ of a curve that can be defined
over $\FF_p$, and so the trace of~\eqref{EQ:basicec} over an even-degree
extension $\FF_q$ of $\FF_p$ is $\pm 2\sqrt{q}$. In particular its trace over
such a field is nonzero. It follows that if we have a curve $D_{u,1}$ over
$\FF_{p^e}$ that is isogenous to its twist, $e$ must be odd, and $u$ must lie in
the intersection of $\FF_{p^e}$ with $\FF_{p^2}$, which is~$\FF_p$.
\end{remark}

\subsection{Doubly isogenous nontwists}
According to the heuristics in Section~\ref{sec:2coverheuristics}, the expected 
number of doubly isogenous nontwists over $\Fq$ should be $\asymp q^{-1/2}$, so
the numbers in the sixth column
of Table~\ref{table:doublyplus} should behave like $2^{-n/2}$. In particular,
as we move from one row to the next, the values in these columns should decrease
by a factor of $\sqrt{2}$. This does not seem to be reflected in the data; we
seem to have more doubly isogenous nontwists than expected. Likewise, we would
expect the numbers in the seventh column of Table~\ref{table:doublyplus} to 
decrease by a factor of $2^{3/2}$ from one row to the next, and this also does
not appear to be happening.

One possible explanation, paralleling the behavior seen in \cite{ArulEtAl2024}, 
is that there may be relations on the parameters $r_1,r_2$ which cause isogenies
between at least three of the elliptic curves appearing in the isogeny
decompositions of the Jacobians of $C_{r_1,1}^{[2]}$ and $C_{r_2,1}^{[2]}$. It
turns out that this does not happen in our case (see 
Theorem~\ref{theorem:finitely many doubly isogenous}), but surprisingly 
Heuristic~\ref{H:default} is still incorrect about doubly isogenous nontwists.
The actual explanation, given by Theorem \ref{T:extraordinary} below, is the
existence of the ``extraordinary curves'' mentioned in the introduction:
two curves $C_1 \colonequals C_{r_1,1}$ and $C_2 \colonequals C_{r_2,1}$ with
$r_1,r_2 \in \bar\QQ$, such that $C_1$ and $C_2$ are doubly isogenous.

\begin{theorem}
\label{T:extraordinary}
Let $K = \QQ(\sqrt{29})$ and let $r_1$ and $r_2$ be the roots of $x^2 - 27x + 1$
in~$K$. Let $C_1 \colonequals C_{r_1,1}$ and $C_2 \colonequals C_{r_2,1}$, where
we use Notation~\textup{\ref{N:basiccurve}}. Let $L$ be the $S_3$
extension of $K$ obtained by adjoining the roots of $x^3 + x^2 + 2$ and 
$x^2 + 1$. Then the Weierstrass points of $C_1$ and $C_2$ are rational over~$L$,
and $C_1$ and $C_2$ are doubly isogenous over $L$ when we embed them into their 
Jacobians using a Weierstrass point as the base point. Furthermore, the Pryms of
the intrinsic triple covers of $C_1$ and $C_2$ are also isogenous over~$L$.
\end{theorem}

\begin{remark}
The field $L$ is the Hilbert class field of the imaginary quadratic field
$\QQ(\sqrt{-29})$. In addition to being the splitting field over $\QQ$ of the
set of polynomials $\{x^2-29,x^2+1,x^3 + x^2 +2\}$, the field $L$ can also be
obtained as the splitting field over $\QQ$ of the single polynomial
$x^{12} + 5 x^{10} - 6 x^8 + 29 x^6 - 6 x^4 + 5 x^2 + 1.$
\end{remark}

Since $L/\QQ$ is a Galois extension of degree~$12$, the natural density of
the rational primes that split completely in $L$ is $1/12$.
For each such splitting prime $p$,
the reductions of $C_1$ and $C_2$ modulo $p$ will have all of their
Weierstrass points rational over $\FF_p$ and will be doubly isogenous. Thus, the
reductions of the curves $C_1$ and $C_2$ will show up as a doubly isogenous pair 
over~$\FF_p$. In addition, the reductions will also have isogenous Pryms from
the intrinsic triple covers.

The proof of the theorem specifies isogenies between the elliptic curves in
Table~\ref{table:lambdas} for $C_1$ and those for~$C_2$. It also specifies 
isogenies between the elliptic curves appearing in the Pryms of the intrinsic
triple covers of $C_1$ and elliptic curves in the corresponding Pryms for~$C_2$.
By checking how twisting $C_1$ and $C_2$ affects these elliptic curves, we see 
that the twists of the reductions of $C_1$ and $C_2$ modulo $p$ will also be
doubly isogenous. Likewise, if $p\equiv 1 \bmod 3$, the intrinsic triple covers of
the twists of the reductions will give rise to isogenous Pryms. However, if 
$p\equiv 2\bmod 3$, the Pryms of the intrinsic triple covers of the twists of the
reductions will not be isogenous, unless one particular elliptic factor of the
Prym has trace~$0$.

Therefore, out of the $1024$ primes closest to $2^n$, we expect that 
$1024/12 = 85\frac{1}{3}$ will have two doubly isogenous pairs coming from the 
reductions of the extraordinary curves of Theorem~\ref{T:extraordinary}. Thus,
the expected number of doubly isogenous pairs in a given row of 
Table~\ref{table:doublyplus} coming from the extraordinary curves is 
$170\frac{2}{3}.$ When we count how many of these doubly isogenous pairs also
have isogenous Pryms from the intrinsic triple covers, we expect to find 
$(2/24)\cdot 1024$ from the primes that are $1\bmod 3$, and $(1/24)\cdot 1024$
from the primes that are $2\bmod 3$, for a total of $1024/8 = 128$. 

Before we give the proof of Theorem~\ref{T:extraordinary}, let us see what our
data shows when we remove the contributions of the extraordinary pair. In 
Table~\ref{table:noextraordinary} we reproduce the data in 
Table~\ref{table:doublyplus} for doubly isogenous nontwists and for doubly
isogenous nontwists with isogenous Pryms coming from the intrinsic triple covers.
We also include columns showing how many pairs of each type are reductions of 
the extraordinary pair (or twists of these reductions), and how many are not.

\begin{table}[t]
\centering
\caption{Data for doubly isogenous curves that are not twists of one another.
For each $n$, the second column lists the total number of (unordered) pairs of 
doubly isogenous curves over $\Fq$ that are not twists of one another, for the 
$1024$ primes $q$ closest to $2^n$. The third column lists the number of these
examples that are (twists of) the reductions of the extraordinary curves from 
Theorem~\ref{T:extraordinary}, and so are explained by the theorem. The fourth
column lists the number that are not explained by Theorem~\ref{T:extraordinary}.
The sixth column lists the number of doubly isogenous pairs from the second 
column that also have isogenous Pryms coming from the intrinsic triple covers, and
the seventh and eighth columns again show how many of these examples are
explained by the extraordinary pair and how many are not.}
\label{table:noextraordinary}
\begin{tabular}{cccrcccrc}
\toprule
&& & \multicolumn{2}{c}{Theorem~\ref{T:extraordinary}?}&& & \multicolumn{2}{c}{Theorem~\ref{T:extraordinary}?}\\
\cmidrule(lr){4-5}\cmidrule(lr){8-9}
  $n$ &&    $2$ & \ \ Yes   & No     && $3$a  & \ \ Yes   & No    \\
\midrule
 $12$ &&     $2568$ & $156$ &           $2412$ && $150$ & $118$ &     $32$  \\
 $13$ &&     $2294$ & $174$ &           $2120$ && $145$ & $128$ &     $17$  \\
 $14$ &&     $1749$ & $176$ &           $1573$ && $135$ & $130$ & \pz  $5$  \\
 $15$ &&     $1428$ & $176$ &           $1252$ && $135$ & $135$ & \pz  $0$  \\
 $16$ &&     $1064$ & $176$ & \pz        $888$ && $133$ & $133$ & \pz  $0$  \\
 $17$ && \pz  $860$ & $176$ & \pz        $684$ && $133$ & $132$ & \pz  $1$  \\
 $18$ && \pz  $631$ & $164$ & \pz        $467$ && $117$ & $117$ & \pz  $0$  \\
 $19$ && \pz  $571$ & $206$ & \pz        $365$ && $151$ & $151$ & \pz  $0$  \\
 $20$ && \pz  $422$ & $178$ & \pz        $244$ && $135$ & $135$ & \pz  $0$  \\
 $21$ && \pz  $328$ & $146$ & \pz        $182$ && $106$ & $106$ & \pz  $0$  \\
 $22$ && \pz  $282$ & $160$ & \pz        $122$ && $114$ & $114$ & \pz  $0$  \\
 $23$ && \pz  $270$ & $172$ & \pz\pz      $98$ && $129$ & $129$ & \pz  $0$  \\
 $24$ && \pz  $236$ & $174$ & \pz\pz      $62$ && $131$ & $131$ & \pz  $0$  \\
 $25$ && \pz  $222$ & $172$ & \pz\pz      $50$ && $132$ & $132$ & \pz  $0$  \\
 $26$ && \pz  $232$ & $204$ & \pz\pz      $28$ && $144$ & $144$ & \pz  $0$  \\
 $27$ && \pz  $194$ & $168$ & \pz\pz      $26$ && $132$ & $132$ & \pz  $0$  \\
 $28$ && \pz  $210$ & $198$ & \pz\pz      $12$ && $153$ & $153$ & \pz  $0$  \\
 $29$ && \pz  $150$ & $146$ & \pz\pz\pz    $4$ && $107$ & $107$ & \pz  $0$  \\
 $30$ && \pz  $183$ & $174$ & \pz\pz\pz    $9$ && $124$ & $124$ & \pz  $0$  \\
\bottomrule
\end{tabular}
\end{table}

\begin{remark}
When we look at the numbers of pairs of doubly isogenous nontwists that are
\emph{not} explained by the extraordinary pair in characteristic~$0$, we see
that this data better fits our heuristic expectations. First, the entries in the
fourth column are plausibly decreasing by a factor of $\sqrt{2}$ from one row to
the next, and second, the entries in the seventh column do rapidly go to~$0$, as
the heuristic says they should.
\end{remark}

\begin{proof}[Proof of Theorem~\textup{\ref{T:extraordinary}}]
For $i=1$ and $i=2$ let $E_i$ be the curve $E_{r_i,1}$ from 
Corollary~\ref{P:baseE}. We check that
\[
x^2 + \frac{4r_1 - 74}{5} x + \frac{1251-35 r_1}{25}
\]
is a factor of the $5$-division polynomial for $E_1$, corresponding to the 
$x$-coordinates of the nonzero points of a $K$-rational subgroup $G$ of $E_1$ of
order~$5$. We also check (by computer algebra, or by Velu's formulas) that the
quotient of $E_1$ by $G$ is isomorphic to $E_2$. Therefore, the Jacobians of 
$C_1$ and $C_2$ are isogenous to one another over~$K$. 

Now let $u_1\in L$ be the $x$-coordinate of one of the Weierstrass points 
of~$C_1$, so that $u_1$ is a root of
\[
x^{12} - 9 x^{10} - 53 x^8 + 266 x^6 + 1707 x^4 + 2183 x^2 + 1.
\]
The Galois group of $L/\QQ$ is isomorphic to~$D_6$; let $\rho$ be the
nontrivial element of the center of this Galois group. (The automorphism 
$\rho$ fixes the roots of $x^3 + x^2 + 2$ and $x^2 + 29$, but negates the 
roots of $x^2 + 1$ and $x^2 - 29$.) If we set $u_2 \colonequals \rho(u_1)$, 
then $u_2$ is the $x$-coordinate of a Weierstrass point of~$C_2$.

Let us view each $C_i$ as being $D_{u_i,1}$, using Notation~\ref{N:ucurves}.
Consider the Orbit~6 curves for $C_1$. If we set 
\[
\lambda_1 \colonequals \frac{4u_1}{(u_1-1)(u_1+3)}
\]
as in Table~\ref{table:lambdas}, then each Orbit~6 curve is a twist of the
elliptic curve
\[
F_1\colon y^2 = x(x-1)(x-\lambda_1)
\]
by the values of $e$ given in the table. Let us denote the $e$-twist of $F_1$ 
by~$F_1^{(e)}$, and let $F_2$ be the image of $F_1$ under~$\rho$.

The $7$-division polynomial for $F_1$ has a unique cubic factor over~$L$, which
corresponds to an $L$-rational subgroup $G$ of $F_1$ of order~$7$. We check that
for each value of~$e$, the quotient of $F_1^{(e)}$ by this subgroup is 
isomorphic to $F_2^{(\rho(e))}$. Thus, the multiset of traces of the elliptic
curves in Orbit~6 for $C_1$ is equal to the corresponding multiset for~$C_2$.

We check that the same holds for Orbit~\threea, Orbit~\threeb, and Orbit~\threec: Each curve in
the given orbit for $C_1$ is $7$-isogenous to the corresponding curve for~$C_2$.

Therefore, $C_1$ and $C_2$ are doubly isogenous over~$L$.

Next, we turn to the intrinsic triple covers of $C_1$ and $C_2$, as in 
Section~\ref{ssec:specialtriple}. Lemma~\ref{L:specialPrym} says that the Prym of 
the intrinsic triple cover of $C_1^{\phihat,1}\to C_1$ is isogenous to the product of the
elliptic curve
\[
F_{\onea}\colon y^2 = x^3 + 81 x^2 + 72 r_1 x + 16 r_1^2
\]
with its twist by $81 - 3r_1$, and since $81-3r_1 = 3(r_1-26)^2/25$ we see that the Prym
is isogenous to the product of $F_{\onea}$ with its twist by~$3$. Similarly, 
Corollary~\ref{C:specialPrym} says that the Prym of the intrinsic triple cover 
$C_1^{\varphi,2}\to C_1$ is isogenous to the product of the elliptic curve
\[
F_{\oneb}\colon y^2 = x^3 + 81 x^2 + 72 (729/r_1) x + 16 (729/r_1)^2
\]
with its twist by $27/r_1- 1$. Since $27/r_1 - 1 = (27-r_1)^2$, this Prym is 
isogenous to the square of~$F_{\twoa}$.

Likewise, the Prym of one of the intrinsic triple covers of $C_2$ is isogenous to
the product of a curve $F_{\twoa}$ with its twist by~$3$, while the Prym of the
other is isogenous to the square of a curve $F_{\twob}$.

As above, we check that $F_{\onea}$ is $3$-isogenous over $L$ to the twist of 
$F_{\twoa}$ by~$3$, so the product of $F_{\onea}$ with its twist by $3$ is isogenous 
to the product of $F_{\twoa}$ with its twist by~$3$. Also, $F_{\oneb}$ is 
$5$-isogenous over $L$ to $F_{\twob}$. Thus, the Pryms of the two intrinsic triple
covers of $C_1$ are isogenous to the Pryms of the two intrinsic triple covers 
of~$C_2$.
\end{proof}

\begin{remark}
Magma code verifying some of the calculations in the above proof are available
in the GitHub repo for this paper~\cite{DoublyIsogenousRepo}.
\end{remark}

\section{An unlikely intersection in characteristic zero}
\label{sec:unlikely}

We begin by recalling the notion of an isogeny correspondence and some of the terminology needed to state the Zilber--Pink conjecture.  

\subsection{Isogeny correspondences} 
\label{ss:isogenycorrespondences}
We work over a field $K$ of characteristic zero. 

\begin{definition}
Let $E^{(1)}$ and $E^{(2)}$ be elliptic curves over $\bA^1$.
An isogeny correspondence between $E^{(1)}$ and $E^{(2)}$ is an irreducible closed curve $Z \subset \bA^1 \times \bA^1$ such that $\pi_1^*(E^{(1)})$ and $\pi_2^*(E^{(2)})$  are isogenous as relative elliptic curves over $Z$, where for $i=1,2$ the $\pi_i$ are the projections $Z \to \bA^1 \times \bA^1$ in the $i$-th coordinate.
In particular, $E^{(1)}_s$ is isogenous to $E^{(2)}_t$ for all $z = (s,t) \in Z$.
\end{definition}

\begin{remark}
This is less general than the definition considered by Buium \cite{Buium1993}.  It suffices for our purposes.  It is also harmless to remove additional points from $\bA^1$ (and also from $Z$) over which the families are not defined.
\end{remark}

\begin{definition} \label{D:simultaneousisogeny}
Let $E^{(1)}, E^{(2)}, E^{(3)}, E^{(4)}$ be elliptic curves over $\bA^1$.  A simultaneous isogeny correspondence between $(E^{(1)}, E^{(2)})$ and $(E^{(3)},E^{(4)})$ is an irreducible closed curve $Z \subset \bA^1 \times \bA^1$ that gives an isogeny correspondence between $E^{(1)}$ and $E^{(3)}$ and also between $E^{(2)}$ and $E^{(4)}$.
\end{definition}

Using modular polynomials, it is straightforward to search for isogeny correspondences between elliptic curves over $\bA^1$ where the isogeny has a fixed degree. For example, suppose $K$ is algebraically closed and we are looking for an isogeny correspondence of degree $N$ between the elliptic curves $E_{s,1}$ and $E_{t,1}$ occurring in Corollary~\ref{P:baseE} for the $D_6$-family.  Plugging the $j$-invariants into the modular polynomial $\Phi_N(x,y)=0$ (\cite[\S69]{Weber}) cuts out a curve (possibly singular or with multiple components) in the $(s,t)$-plane that gives a degree-$N$ isogeny correspondence between $E_{s,1}$ and $E_{t,1}$. 

Carrying out these computations for small $N$ is the approach used in \cite[\S6]{ArulEtAl2024} to find a variety of simultaneous isogeny correspondences for the elliptic curves showing up in the Jacobian of the double cover of the $D_4$-family.  We do not find interesting simultaneous isogeny correspondences for the $D_6$-family.  However, ruling them out requires additional ideas as we cannot use this method for infinitely many $N$; this is the subject of Section~\ref{sec:computingisogeny}.

\begin{remark}
\label{remark:simultaneous}
Consider looking for doubly isogenous curves $C_1$ and $C_2$ in the \dsixfamily. This requires at least five simultaneous isogenies between elliptic curves:
between the elliptic factor of $\Jac C_1$ and that of $\Jac C_2$,
and between the four elliptic curves appearing in the Prym of $C_1^{[2]}\to C_1$ and the four appearing in the Prym of $C_2^{[2]}\to C_2$; see 
Example~\ref{example:doubleisogenypatterns}.
Each correspondence is a plane curve, so we would expect that imposing two isogeny conditions would limit us to points $(s,t)$ in the intersection of two plane curves.  We expect this intersection to be a finite set.  Furthermore, we would expect that three or more isogeny conditions will leave us with no points $(s,t)$, as the intersection of three generic plane curves is empty.
From this viewpoint, the extraordinary curves in Theorem~\ref{T:extraordinary} represent a point in the intersection of \emph{five} plane curves!
\end{remark}

One complication is that we may vary the degree of the isogeny, so there are multiple possible curves which realize an isogeny.  This naturally leads us to the Zilber--Pink conjecture.

\subsection{The Zilber--Pink conjecture}
We now work over an algebraically closed field $K$ of characteristic zero. We state the conjecture in its general form but omit the definitions. See \cite{BarroeroDill} for a general discussion and Example \ref{ex:ZP} for the precise context and definitions that we will require.

\begin{conjecture}
\label{conj:ZP}
Let $V$ be an irreducible subvariety of a mixed Shimura variety~$S$. If $V$ is
not contained in a special subvariety of positive codimension, then the intersection
of $V$ with the union of all special subvarieties of $S$ of codimension greater than $\dim V$ is not Zariski dense in $V$.
\end{conjecture}

This remains open, but some progress has been made in the case that $V$ is a curve. Unfortunately, we need apply it to a surface, and this case appears to be still open (J. Pila, personal communication). See \cite{DawOrr2025,HabeggerPila2012} for background and a review of special subvarieties when $S$ is a product of modular curves, which is the relevant case for us. These two papers also prove results for the case that $V$ is a curve, which may suffice for the second invocation of the conjecture in the proof of Theorem \ref{theorem:finitely many doubly isogenous}; we have not checked this, though, because the first invocation of the conjecture in the proof is for a surface. 

\begin{example}
\label{ex:ZP}
In the case of interest, we work over $K = \overline{\QQ}$ and take $S = Y(1)^n$ for some $n$, where $Y(1)$ is the affine line, viewed as the coarse moduli space of elliptic curves.  A \emph{special} subvariety of $Y(1)^n$ is a subvariety defined by equations of the form
\begin{itemize}
\item  $\Phi_N(j_a,j_b) =0$ for some integer $N \geq 1$ and $1 \leq a < b \leq N$; or
\item  $j_a = c$ for some $1 \leq a \leq N$ and CM $j$-invariant $c \in \overline{\QQ}$.
\end{itemize}
\end{example}

\begin{theorem} \label{theorem:finitely many doubly isogenous}
Conjecture~\textup{\ref{conj:ZP}} in the context of Example~\textup{\ref{ex:ZP}} implies that, in characteristic zero, there are only finitely many pairs of doubly isogenous curves in
the \dsixfamily \eqref{EQ:basiccurve}.
\end{theorem}

\begin{proof}
Two curves in the \dsixfamily are doubly isogenous if the elliptic curves in Corollary~\ref{P:baseE} corresponding to the base curve are isogenous and if there are isogenies between the elliptic curves that occur when decomposing the Prym up to isogeny.   Recall that the pattern of isogenies may be more complicated than simply requiring isogenies between the corresponding orbits (see Example~\ref{example:doubleisogenypatterns}).  For any such pattern, we may select elliptic curves $E^{(1)}, E^{(2)}, \ldots, E^{(6)}$ depending nontrivially on a parameter such that 
 a pair of doubly isogenous $D_6$-curves for parameters $t$ and $s$ gives us two triples of elliptic curves
$$(E^{(1)}_t,E^{(2)}_t,E^{(3)}_t),(E^{(4)}_s,E^{(5)}_s,E^{(6)}_s)$$ 
with $E^i_t \sim E^{3+i}_s$ for $i=1,2,3$.  (See Table \ref{table:lambdas} for the expression of the $\lambda$-invariant for each elliptic curve and the text immediately following for the corresponding $j$-invariants, and see Remark \ref{remark:simultaneous} for an explanation of the minimal choice of $6$ curves.)   One example would be taking $E^{(1)} = E^{(4)}$ to be the base elliptic curve and the remaining elliptic curves to be from two of the orbits that need to be isogenous.

 The triple of $j$-invariants for $(E^{(1)},E^{(2)},E^{(3)})$ gives a curve in $Y(1)^3$, and similarly $(E^{(4)},E^{(5)},E^{(6)})$ gives a curve in $Y(1)^3$.  By taking the closure,   we obtain irreducible closed curves $X, X' \subset Y(1)^3$.
Consider the surface $X \times X' \subset Y(1)^6$.  
To establish the theorem, it suffices to show there are only finitely many points $(j_1,j_2,\ldots,j_6)$ on $X \times X'$ such that there exist positive integers $n_1,n_2,n_3$ such that $\Phi_{n_i}(j_{i},j_{i+3})=0$ for $i=1,2,3$ (equivalently, there exists an isogeny of degree $n_i$ between the corresponding elliptic curves).  

Notice that $X \times X'$ is not a special subvariety:
\begin{itemize}
\item  $X \times X'$ is not contained in any hyperplane of the form $j_i = a$ since the $E^i$ depend nontrivially on the parameter.

\item Remark~\ref{R:noisogenies} shows that $X \times X'$ is not contained in any hypersurface of the form $\Phi_N(j_a,j_b)=0$ since there exist choices of the parameters in the $D_6$-family where all ten elliptic curves are non-isogenous.
\end{itemize}
Then the Zilber--Pink conjectures implies that $(X \times X') \cap S^{[3]}$ is not Zariski dense in $X \times X'$, where $S^{[m]}$ is the union of all special subvarieties of $S$ of codimension at least $m$.
The three conditions $\Phi_{n_i}(j_{i},j_{i+3})=0$ define a special subvariety in $Y(1)^6$ of dimension three.  Thus the doubly isogenous curves in the $D_6$-family give points on $(X \times X') \cap S^{[3]}$.  It suffices to show there are no irreducible one-dimensional components $V$ in the Zariski closure of $(X \times X') \cap S^{[3]}$ in $X \times X'$.  

For such a $V$, the intersection $V \cap S^{[3]}$ is Zariski dense in $V$, and $V \cap S^{[3]}$ is contained in $V \cap S^{[2]}$, so we can apply
the Zilber--Pink conjecture to deduce that $V$ is special. In particular, as $X \times X'$ is not contained in any hyperplane of the form $j_i=a$, we conclude that $V$ is cut out by (at least) five equations of the form $\Phi_N(j_a,j_b)=0$ for varying positive integers $N$ and $1 \leq a < b \leq 6$.  
If there were two equations $\Phi_{N_1}(j_a,j_b)=0$ and $\Phi_{N_2}(j_a,j_b)=0$ with the same indices, there would be a non-integer endomorphism of $E^a$ and hence $E^a$ would be CM.  This is impossible as the family $E^a$ is non-constant, so there are (at least) five isogenies $E^a \sim E^b$ for different indices.  Thus pulling $V$ back to the $s,t$-plane would give a simultaneous isogeny correspondence between $E^{(1)},E^{(2)},E^{(3)}$ and $E^{(4)},E^{(5)},E^{(6)}$.  These do not exist, which we can verify using the computational techniques in Section~\ref{sec:computingisogeny} (see Example~\ref{ex:nocorrespondenced6}).  
Thus $X \times X' \cap S^{[3]}$ is zero-dimensional, as desired.
\end{proof}
\section{Simultaneous isogeny correspondences}
\label{sec:computingisogeny}

In this section we describe an algorithm for finding all simultaneous isogeny correspondences between families of elliptic curves.  
More precisely, working over a number field $k$, let $E^{(1)},E^{(2)}, E^{(3)}, E^{(4)}$ be elliptic curves over $\bA^1$.  The algorithm will output all simultaneous isogeny correspondences between $(E^{(1)},E^{(2)})$ and $(E^{(3)},E^{(4)})$ (recall Definition~\ref{D:simultaneousisogeny}).  
We have implemented the algorithm in Magma and put our code in a GitHub repository \cite{DoublyIsogenousRepo}.

\subsection{Detecting isogeny correspondences}

As discussed in Section~\ref{ss:isogenycorrespondences}, it is straightforward to search for isogenies of degree $n$ between two families of elliptic curves using the $n$th modular polynomial.   This allows us to computationally find all isogeny correspondences where the degree of the isogenies lie in any fixed finite set of positive integers.  However, it provides no way to rule out the existence of additional isogeny correspondences whose degrees are outside this set.

We make use of work of Buium to detect isogeny correspondences without restriction on the degree.
 Fix a differentially closed field $K$ of characteristic zero with derivation denoted by $y \mapsto y'$.

\begin{definition}
For $ y \in K$, define 
\begin{equation} \label{def:chi}
\chi(y) \colonequals \frac{ 2y' y''' - 3 (y'')^2}{4 (y')^2} + (y')^2 \frac{ y^2 - y +1 }{4 y^2 (y-1)^2}.
\end{equation}
\end{definition}

\begin{theorem}[{Buium \cite[Chapter 6 Theorem 5.2]{Buium1994}}] \label{theorem:buium}
Let $E^{(1)}$ and $E^{(2)}$ be elliptic curves over $K$ with $\lambda$-invariants $\lambda_1, \lambda_2 \in K\setminus \{0,1\}$.  In other words, in Legendre form, $E^{(1)}$ and $E^{(2)}$  are
\[
E^{(1)}\colon  y^2 = x ( x-1)(x-\lambda_1), \quad \text{and} \quad E^{(2)}\colon y^2 = x (x-1)(x-\lambda_2).
\]
If $\chi(\lambda_1) \neq \chi(\lambda_2)$ then $E^{(1)}$ and $E^{(2)}$ are not isogenous over $K$.
\end{theorem}

Now working over the number field $k$, for $i=1,2$ let $\pi_i\colon \bA^1_s \times \bA^1_t \to \bA^1_x$ denote the projection maps where the subscripts label the parameters.  
Let $E^{(1)}$ and $E^{(2)}$ be elliptic curves over $\bA^1_x$ with $\lambda$-invariants $\lambda_1,\lambda_2 \in k(x)$. 

\begin{corollary} \label{corollary:differentialisogeny}
Suppose $\lambda_1,\lambda_2 \neq 0,1$.  If there is an isogeny correspondence ${f\colon Z \subset \bA^1_s \times \bA^1_t}$ between $E^{(1)}$ and $E^{(2)}$, then over $Z$ we have
\begin{equation}
\chi(\lambda_1(s))) = \chi(\lambda_2(t)).
\end{equation}
\end{corollary}

\begin{proof}
There is a natural map $k[s,t] \to k(Z)$, and we take $K$ to be the differential closure of $k(Z)$.  Applying Theorem~\ref{theorem:buium} to $ (\pi_1 \circ f)^*(E^{(1)})$ and $(\pi_2 \circ f)^*(E^{(2)})$ gives the result.
\end{proof}

\subsection{An algorithm for finding simultaneous isogeny correspondences}
\label{sec:isogenyalgorithm}
Given elliptic curves $E^{(1)},E^{(2)}, E^{(3)}, E^{(4)}$ over $\bA^1$ with $\lambda$-invariants $\lambda_1,\lambda_2,\lambda_3,\lambda_4 \in k(\bA^1_k)$, we now turn to producing a complete list of simultaneous isogeny correspondences between $(E^{(1)},E^{(2)})$ and $(E^{(3)},E^{(4)})$ making use of Corollary~\ref{corollary:differentialisogeny}.  
We do so by finding simultaneous solutions to
\begin{equation} \label{eq:simulatenousdiffeq}
\chi(\lambda_1(s)) = \chi(\lambda_3(t)), \quad 
\chi(\lambda_2(s)) = \chi(\lambda_4(t)) ,
\end{equation}
viewing $s$ as a function of $t$.  More formally, we pull back $s$ and $t$ to the (unknown) irreducible curve $Z$ realizing the isogeny correspondence and solve \eqref{eq:simulatenousdiffeq} there.  These are third order differential equations, and the first step is to combine them to reduce the order.  

\begin{lemma}
The equation 
\begin{equation} \label{eq:firstordereq}
\chi(\lambda_1(s)) - \chi(\lambda_3(t)) - \chi(\lambda_2(s))  + \chi(\lambda_4(t))=0
\end{equation}
 is a first order differential equation, which we denote by $F(t,s,s')=0$.
\end{lemma}

\begin{proof}
Let $f = \lambda_i$.  Using the chain rule to compute derivatives of $f(s)$ in terms of $t$, we compute that
\begin{align*}
\chi(f(s)) &= \frac{1}{2 f'(s) s'} \left( f'(s) s''' + 3 f''(s) s' s'' + \cdots \right) \\
& \qquad\qquad - \frac{3}{4} \left( \frac{ \left( f'(s) s''  + f''(s) (s')^2 \right)^2}{ f'(s)^2 (s')^2} \right) + \cdots \\
&= \frac{s'''}{2s'} - \frac{3}{4} \frac{(s'')^2}{(s')^2} + \cdots
\end{align*}
where the omitted terms are rational functions of $s'$, $s$, and $t$.  Note the non-omitted terms are independent of $f$. Furthermore, observe that $\chi(f(t))$ is a rational function of $t$.  Thus the second and third order terms in \eqref{eq:firstordereq} cancel, giving a first order equation.  
\end{proof}

The second step is to differentiate $F(t,s,s')=0$ with respect to $t$, obtaining
\[
0 = \frac{dF}{dt} = \frac{\partial F}{\partial t} + \frac{\partial F}{\partial s} s' + \frac{\partial F}{\partial s'} s''
\]
and hence a formula for $s''$ in terms of $t,s,s'$.  Differentiating again gives a formula for $s'''$.  Substituting them into the equation $\chi(\lambda_1(s)) = \chi(\lambda_3(t))$ gives a first order equation $G(t,s,s')=0$.  

The third step is to compute the resultant of $F$ and $G$ with respect to $s'$.  This is an algebraic equation $R(s,t)=0$ which is satisfied when there is a choice of $s'$ satisfying both $F$ and $G$ (and hence \eqref{eq:simulatenousdiffeq}).  We differentiate $R$ to obtain a formula for $s'$ in terms of $t,s$ and substitute into $F$ to obtain an additional algebraic equations $H(s,t) =0$. 

If there is an isogeny correspondence $Z \subset \bA^1_s \times \bA^1_t$ then $R$ and $H$ will be zero on $Z$ by \eqref{eq:simulatenousdiffeq}.  Therefore the last step is to find common factors of (the numerators of) $R$ and $H$: Each irreducible common factor gives a potential simultaneous isogeny correspondence.  

A common factor does not necessarily give a simultaneous isogeny correspondence.  It is also possible that the algebraic manipulations that produce $R$ and $H$ introduce additional factors not corresponding to isogeny correspondences, so the last step of the algorithm is to use the factor to compute derivatives and check whether \eqref{eq:simulatenousdiffeq} is satisfied.

\begin{example} \label{ex:663a3a}
We can look for simultaneous isogeny correspondences in the elliptic curves occurring in Table~\ref{table:lambdas}.  For example, let us look for simultaneous isogeny correspondences where Orbit~6 is isogenous to Orbit~6, and Orbit~3a isogenous to Orbit~3a (we don't ask for any relationship between the elliptic curves $E_{3b}$ and $E_{3c}$ in other orbits or the elliptic curve $E_{r,c}$ from Corollary \ref{P:baseE}).  
After a computation, we find the common factors of $R$ and $H$ consist of:
\begin{itemize}
\item  Four expected factors involving $s$ and $t$: $s-t$, $s+t$, $st-3$, $st +3$.  Over each zero set, $\pi_1^*(E_6) \sim \pi_2^*(E_6)$ and $\pi_1^*(E_{3a}) \sim \pi_2^*(E_{3a})$.

\item  An irreducible polynomial in $\QQ[s,t]$ of degree $28$ in $s$ and $t$ which does not give an isogeny correspondence.  In particular, using this polynomial to implicitly compute derivatives of $s$ with respect to $t$ gives results that do not satisfy \eqref{eq:simulatenousdiffeq}.  Alternately, specializing $s$ and $t$ subject to this relation and using those values to construct elliptic curves with the specified $\lambda$-invariants gives non-isogenous elliptic curves.  

\item  Factors only involving one variable.  As there are no isogeny correspondences with one variable held fixed, these are automatically spurious.
\end{itemize}

The four expected factors correspond to degree one isogenies which exist for trivial reasons.  In light of Lemma~\ref{lemma:isomorphic D_{u,c}}, all of the curves $D_{u,c}$ with $u \in \{s,-s, 3/s,-3/s\}$ are geometrically isomorphic.  Hence the elliptic curves appearing in Orbit~6 must be geometrically isomorphic as well, as that orbit is preserved by the $D_6$-action. 
Furthermore, these four values of $u$ produce the same value of $16u^2/(u^2+3)^2$, making the elliptic curves in Orbit~3a geometrically isomorphic as well.  This explains the four factors giving expected simultaneous correspondences between $(E_6,E_{3a})$ and $(E_6,E_{3a})$.
\end{example}

\begin{remark}
The simultaneous isogeny correspondences found in Example~\ref{ex:663a3a} actually come from situations where the genus two \dsixcurves $D_{s,c}$ and $D_{t,c}$ are geometrically isomorphic and hence doubly isogenous.  This is not necessarily true for other isogeny correspondences.  
\end{remark}

\subsection{Practical considerations}
While conceptually simple, carrying out this approach is a bit complicated as the rational functions appearing in these algebraic differential equations are quite large, even for relatively simple choices of $\lambda$-invariants~$\lambda_i$.  We also want to consider many instances of this problem in Examples~\ref{ex:nocorrespondenced6} and \ref{ex:correspondenced4}.  Therefore we discuss some computational improvements that make the computations tractable.  

The first simplification is to realize it suffices to compute the numerators of the rational functions $F, G, R, H$ as we are only interested in when they are zero.  The expressions we obtain for the numerator are not guaranteed to be relatively prime to the (uncomputed) denominator, so it is essential to rule out spurious solutions at the end of the computation by computing derivatives and checking the original differential equations.  As illustrated in Example~\ref{ex:663a3a}, this can happen.  The entire computation takes $30$--$60$ minutes on a single core using this approach, and yields an explicit set of potential isogeny correspondences. 

\begin{remark}
 The fundamental obstruction is that (the numerator of) $H$ is legitimately complicated: In Example~\ref{ex:663a3a}, it has degree $1166$ in $s$ and also in $t$, and the coefficients are rational numbers of large height.  The polynomial $R$ is ``simpler,'' having only degree $568$ in each variable.  Note that one variable factors cause the majority of the complexity ($R$~is a multiple of a degree-$504$ polynomial in $\QQ[s]$) but can basically be ignored when looking for simultaneous isogeny correspondences.
 \end{remark}

This simplification is sufficient to search for isogeny correspondences in any single example, and for the computations in Example~\ref{ex:663a3a}.  However, to prove Theorem~\ref{theorem:finitely many doubly isogenous} requires considering many cases, with $\lambda_i$ drawn from a set of size 5 (the $\lambda$-invariants of the elliptic curve from Corollary~\ref{P:baseE} and the 4 orbits in Table~\ref{table:lambdas}).  Even exploiting symmetry, this is a substantial number of cases to consider.   We can further speed up the computations by performing most of the computation after specializing one of the variables.  

\begin{enumerate}
\item  We compute (the numerator of) $R(s,t)$, which is substantially simpler than $H(s,t)$.

\item  We remove a list of known simultaneous isogeny correspondences from the numerator of $R(s,t)$,  as well as all one-variable factors, giving a polynomial $\widetilde{R}(s,t)$. 

\item  We compute $H(s_0,t)$ for an $s_0 \in \QQ$ such that $\deg_t(H(s_0,t)) = \deg_t(H(s,t))$.

\item  Next we compute all common factors $f_{s_0}(t)$ of $H(s_0,t)$ and $\widetilde{R}(s_0,t)$ in $\QQ[t]$.

\item  \label{laststep} For each such factor $f_{s_0}(t) \in \QQ[t]$, we reduce modulo an appropriate prime $p$ and choose a root $t_0$ of $f_{s_0}(t) \equiv 0 \pmod{p}$ over a small extension of $\FF_p$.  We check whether $E_{1,s_0}$ is isogenous to $E_{3,t_0}$ and $E_{2,s_0}$ is isogenous to $E_{4,t_0}$. 

\end{enumerate}

The significance of specializing $H(s,t)$ at a value of $s$ which does not decrease the degree in $t$ is that any common factor of $H(s,t)$ and $\widetilde{R}(s,t)$ with degree $d$ will then specialize to a common factor of $H(s_0,t)$ and $\widetilde{R}(s_0,t)$ of degree~$d$.  This relies on the fact that we already removed all one-variable factors from the numerator of $R(s,t)$.  In particular, any simultaneous isogeny correspondence will show up after specialization.

For step~\eqref{laststep}, it is easy to check for isogenies over finite fields via point counting.  If the elliptic curves in the last step are not isogenous, then we know that $f_{s_0}(t)$ is not the specialization of a simultaneous isogeny correspondence.  If this conclusion holds for all common factors of $H(s_0,t)$ and $\widetilde{R}(s_0,t)$, we conclude there are no additional simultaneous isogeny correspondences beyond those already known.

\begin{remark}
If the elliptic curves in step \eqref{laststep} are isogenous, with high probability $f_{s_0}(t)$ is the specialization of a common factor of $H(s,t)$ and $\widetilde{R}(s,t)$.  However, there is a small chance the isogenies exist after reducing modulo $p$ by coincidence, so we cannot prove the existence of an isogeny correspondence after specialization.  This is why we first remove known isogeny correspondences (found using modular polynomials).  In practice, for consistency checking and to further reduce the probability of false positives we specialize in each variable and perform \eqref{laststep} twice.
\end{remark}

\begin{example}
With these speedups, verifying that the four trivial simultaneous isogeny correspondences identified in Example~\ref{ex:663a3a} are the only ones takes about $30$ seconds.  
\end{example}

\begin{example} \label{ex:nocorrespondenced6}
For the $D_6$-family, it takes around $16$ hours\footnote{Running a single instance of Magma on an Intel Core i9-10900X CPU at 3.70GHz.} to rule out the existence of a nontrivial simultaneous isogeny correspondence among any combination of pairs of the five elliptic curves occurring when decomposing the Jacobian of~$D^{[2]}_{u,c}$.   This computation is the remaining ingredient in the proof of Theorem~\ref{theorem:finitely many doubly isogenous}.  The computation in Example~\ref{ex:663a3a} is one of the $200$ cases to be considered.  Another example case is looking for a simultaneous isogeny correspondence between $(E_6,E_{3a})$ and $(E_6,E_{3b})$, and computationally verifying there are no nontrivial ones.   Note that we have already taken advantage of the $D_6$-action to reduce the number of elliptic curves  under consideration from $16$ (the one occurring in the Jacobian of the base and $15$ occurring in the Prym of $D^{[2]}_{u,c} \to D_{u,c}$) to five (the one occurring in the Jacobian of the base and the four from the four orbits in Table~\ref{table:lambdas}).
\end{example}

\begin{example} \label{ex:correspondenced4}
For genus-$2$ curves with $D_4$-action, \cite{ArulEtAl2024} identifies a number of simultaneous isogeny correspondences which affect the heuristics for doubly isogenous curves.  As discussed in \cite[\S6]{ArulEtAl2024}, these were found by searching for low-degree isogenies using modular polynomials. 
It was not clear whether \emph{all} simultaneous isogeny correspondences were found, although as discussed in \cite[\S7]{ArulEtAl2024} the ones found were sufficient to explain the experimental data about doubly isogenous curves. 

Using the techniques of this section one can easily check whether there are additional simultaneous isogeny correspondences.  In fact, the highest degree of isogeny found is $4$, so all simultaneous isogeny correspondences are in the range already searched \cite[Remark 6.5]{ArulEtAl2024}.  
\end{example}

\section{Application to factoring polynomials}
\label{sec:factoring}

Poonen showed how to use a suitable one-parameter family of abelian varieties $\CA$ defined by equations with coefficients in $\ZZ[t]$ to factor polynomials over finite fields in deterministic polynomial time 
\cite[Theorem 5.1]{Poonen2019}.
A one-parameter family is suitable if, for all sufficiently large primes $p$, the specializations $\CA_u$ at generic $u\in \FF_p$ have distinct zeta functions; here \emph{generic} means $\Delta(u)\ne 0$, for some fixed $\Delta\in \ZZ[t]$ defining the bad locus of $\CA$.
An explicit candidate for such an $\CA$ was given in \cite[\S5]{SutherlandVoloch2019} via the one-parameter family of genus-4 hyperelliptic curves
\[
X_t\colon y^2 = x^9 + (t-1)x^3+tx^2+(t+1)x+1,
\]
whose Jacobians at generic $u\in \FF_p$ have distinct zeta functions for $p\le 2^{17}$, and conjecturally for all primes $p$.

Our \dsixfamily of doubly isogenous curves leads to an alternative candidate for $\CA$ that is computationally more tractable (even though of higher dimension), because the Jacobians are isogenous to products of elliptic curves.
Recall from Proposition~\ref{P:2decomp} that there are five geometrically distinct elliptic curve factors of the Jacobian of the unramified elementary abelian $2$-cover $\pi\colon D^{[2]}_{u,c}\to D_{u,c}$.
These are the elliptic curve $E_{r,c}$ from Corollary \ref{P:baseE} and the elliptic curves with $j$-invariants $j_6, j_\threea, j_\threeb, j_\threec$ described in Table~\ref{table:lambdas} in terms of the parameters $u$ and $c$.  By Lemma~\ref{L:rationalW}, the parameters $r$ and $u$ are related via
\[
r = -\frac{u^2(u^2-9)^2}{(u^2-1)^2}.
\]
Fixing $c=1$, these five elliptic curves depend on a single parameter $u$, and for $u\not\in U_\sing \subseteq \{0,\pm 1,\pm 3,\pm\sqrt{-3}\}$ they are all nonsingular (see Notation~\ref{N:ucurves}).

Taking $\CA$ to be the Jacobian of the unramified cover $\pi$ (or more parsimoniously, the product of these five elliptic curves) will not suffice for two reasons.  The first reason is the existence of doubly isogenous \dsixcurves in infinitely many characteristics (Theorem~\ref{T:extraordinary}).
The second reason this does not suffice is that $12$ values of the parameter $u$ will give the same $D_6$-curve (Lemma~\ref{lemma:isomorphic D_{u,c}}): The parameter $u$ corresponds to the choice of a rational Weierstrass point.
Recall that because of the $D_6$-action, all of the maximal unramified abelian $2$-covers (which depend on the choice of a Weierstrass point) are isomorphic.  Thus there is automatically duplication due to the parametrization.  

To overcome this, we consider additional elliptic curves obtained by shifting the original five by $u\mapsto u-22$ and by $u\mapsto u-30$.
This yields a total of $15$ elliptic curves, each of which can be put in the form $y^2=x^3+Ax+B$ with $A,B\in \ZZ[u]$.
We could take $\CA$ to be the product of these 15 elliptic curves, but for the purpose of our factorization algorithm, we regard them as an \emph{ordered} tuple, as we will compute with each elliptic curve factor separately.
We may specialize at $u\in \FF_p$ when 
\begin{equation}
\label{eq:badfactoring}
u,u-22,u-30\not\in U_\sing \subseteq \left\{0,\pm 1,\pm 3,\pm\sqrt{-3}\right\},
\end{equation}
obtaining a $15$-tuple of Frobenius traces (equivalently, zeta functions) for the $15$ elliptic curves. 
(It is a coincidence that $15$ is also the dimension of the Prym variety associated to the double cover of $D_{u,1}$; some of the elliptic curves do not occur in this Prym.)

\begin{conjecture} \label{conj:distincttuples}
For all primes, generic $u \in \FF_p$ yield distinct $15$-tuples of Frobenius traces.
\end{conjecture}

We have verified Conjecture \ref{conj:distincttuples} for all primes $p\le 2^{20}$.
We emphasize that the traces in a single tuple need not be distinct; it is the tuples that are distinct.  

\begin{remark}
We use two shifts as we expect this to remove duplication due to the parametrization.  In particular, given $u$ there are four choices of $u'$ such that the
$j$-invariants of the five fundamental elliptic curves at $u$ and $u'$ are the same.
(There are four, not twelve, as the $D_6$-action permutes orbits \threea, \threeb, and \threec.
Some of the symmetry is broken by ordering the elliptic curves.)  View this as a reducible curve in the $u,u'$ plane.
Do the same for shifts of $u$ by $22$ and $30$.
We do not expect any common points on these three plane curves for dimension reasons.
Common points would correspond to values of $u$ and $u'$ which give the same $15$-tuple of $j$-invariants and hence Frobenius traces.
Thus adding two shifts heuristically removes the duplication due to parametrization.

It turns out that this shifting also appears to suffice to distinguish the doubly isogenous curves.  We chose to shift by $22$ and $30$ as this choice did not result in any coincidences for small primes.
\end{remark}

We now turn to describing the factoring algorithm.  As noted in \cite{Poonen2019}, Poonen's result is inspired by a 2005 suggestion of Kayal in which $\CA$ is an elliptic curve, a setting in which the factorization algorithm is particularly efficient.
The necessary hypothesis that the specializations $\CA_u$ have distinct zeta functions for $u\in \FF_p$ cannot be met when $\CA$ is an elliptic curve (there are not enough zeta functions to go around), but it can when $\CA$ is a product of elliptic curves, and this is nearly as efficient.
While both can be accomplished in deterministic polynomial time, computing the Frobenius traces of 15 elliptic curves  over~$\FF_p$ is much easier than computing the zeta function of a genus-$4$ curve over~$\FF_p$.
The former can be done in $O(\log^5 p)$ time via Schoof's algorithm \cite{Schoof1985} (this complexity bound follows from \cite[Corollary 11]{ShparlinskiSutherland} and \cite{HarveyVanDerHoeven}), while the most efficient methods known for computing the zeta function of a genus-$g$ hyperelliptic curve $C$ take $O(\log^{cg} p)$ time with $c\ge 4$, in general \cite{AbelardGaudrySpaenlehauer}.
Furthermore, we are not aware of any practical implementation of this algorithm for $g>3$.
The algorithm is also easier to describe in the elliptic curve setting, so let us briefly recall Kayal's approach.

As shown by Berlekamp \cite{Berlekamp1970}, the problem of factoring an arbitrary polynomial over an arbitrary finite field $\FF_q$ can be efficiently and deterministically reduced to the problem of factoring a monic polynomial $f$ over a prime field $\FF_p$ that splits completely into a product of distinct linear factors.
It suffices to find any nontrivial factorization of $f$, since we can then recursively factor each piece. 

\begin{remark}
There is a simple probabilistic approach to factoring such an~$f$, due to Rabin \cite{Rabin1980}, that efficiently computes $\gcd(f(v),(v+\delta)^{(p-1)/2}-1)$ for random $\delta\in \FF_p$; this can be done in time that is quasilinear in $d$ and quasiquadratic in $\log p$.
On average only two values of~$\delta$ need to be tested before a nontrivial factor of $f$ is found.
This effectively solves the problem from a practical perspective, but no deterministic algorithm is known to solve this problem in polynomial time, not even in the simplest possible case where $f(v)=(v-a)(v-b)$ has degree $d=2$.
(Under the generalized Riemann hypothesis, the case $d=2$, and indeed any case with $d=O(1)$, can be solved in deterministic polynomial time \cite{Ronyai1992}; but the problem is open in general, even under the generalized Riemann hypothesis.)
\end{remark}

To factor $f$ deterministically, Kayal proposed a modification of Schoof's algorithm for computing the trace of Frobenius of an elliptic curve $E/\FF_p$: Instead of working over $\FF_p$, Kayal works over the \'etale $\FF_p$-algebra $\CF\colonequals \FF_p[v]/(f(v))\simeq \FF_p^d$, and checks for nontrivial zero divisors in $\CF$ using gcds, as we now explain.

We first recall Schoof's algorithm for computing the trace of Frobenius $a_p$ of an elliptic curve $E/\FF_p$, with $p>3$ and $E\colon y^2=x^3+Ax+B$ in short Weierstrass form. Schoof's algorithm computes $a_p$ modulo primes $\ell=2,3,5,\ldots$ with $\ell\ne p$ by testing identities in the $\ell$-torsion endomorphism ring $\End(E[\ell])$ of the form
\begin{equation}\label{eq:schoof}
\pi^2-[c]\pi+[p_\ell]=0.
\end{equation}
Here $\pi$ denotes the Frobenius endomorphism $(x\,{:}\,y\,{:}\,z)\mapsto (x^p\,{:}\,y^p\,{:}\,z^p)$, while $[c]$ and $[p_\ell]$ denote scalar multiplication by $c$ and $p_\ell\colonequals p\bmod \ell$.  This equation will be satisfied if and only if $c\equiv a_p\bmod \ell$, and one can simply check $c=0,1,2,\ldots,\ell-1$ until \eqref{eq:schoof} is satisfied.
Once $a_p\bmod \ell$ is known for sufficiently many $\ell$ (we need the product of the $\ell$'s to exceed the Hasse bound of $4\sqrt{p}$), the value of $a_p\in\ZZ$ is obtained via the Chinese remainder theorem.  As noted above, this yields a deterministic algorithm that computes $\#E(\FF_p)$ in $O((\log p)^{5})$ time.

As explained for example in \cite[Lecture~8]{Sutherland18783}, each nonzero endomorphism in $\End(E[\ell])$ is uniquely represented by a pair $(g(x),h(x)y)$ of elements in the ring $\FF_p[x,y]/(x^3+Ax+B-y^2,\phi_\ell(x,y))$, where $\phi_\ell$ is the $\ell$th division polynomial of $E$ and $g,h\in \FF_p[x]$ are reduced modulo $\phi_\ell$.
For $c=0$ the identity \eqref{eq:schoof} holds if and only if $g_1=g_2$ and $h_1=h_2$, where $g_1,h_1,g_2,h_2\in \FF_p[x]$ are the unique polynomials of degree less than $\deg\phi_\ell$ for which $(g_1(x),h_1(x)y)$ is the unique representation of $\pi^2$ in $\End(E[\ell])$ and $(g_2(x),h_2(x)y)$ is the unique representation of $-[p_\ell]$ in $\End(E[\ell])$; for $0<c<\ell$
one uses $(g_1(x),h_1(x)y)$ representing $\pi^2+[p_\ell]$ and $(g_2(x),h_2(x)y)$ representing $[c]\pi$.
In the brief summary of Schoof's algorithm above we have ignored issues that may arise when $\phi_\ell$ is not irreducible, but \cite[Lecture~8]{Sutherland18783} explains how to efficiently handle this possibility.

When working over the $\FF_p$-\'etale algebra $\CF\colonequals \FF_p[v]/(f(v))\simeq \FF_p^d$, rather than testing $g_1=g_2$ and $h_1=h_2$ one instead computes $\gcd(f,g_1-g_2)$ and $\gcd(f,h_1-h_2)$.
If \eqref{eq:schoof} holds in some but not all factors of $\CF$, this reveals a proper divisor of~$f$.
This is guaranteed to happen if the reductions of $E\colon y^2 = x^3 + A(u)x + B(u)$ to two distinct factors of $\CF$ have different Frobenius traces.
The complexity of running Schoof's algorithm over $\CF$ is $O(d^{1+o(1)}(\log p)^{5+o(1)})$, which is polynomial in the size of~$f$, which is $O(d\log p)$ bits.

As noted above, this approach cannot be made to work using a single elliptic curve  $E\colon y^2 = x^3 + A(u)x + B(u)$, but under Conjecture \ref{conj:distincttuples} it is guaranteed to work if we run the same algorithm on fifteen elliptic curves (which amounts to computing the zeta function of their product).
To factor the polynomial $f\in \FF_p[v]$ defining the \'etale algebra $\CF=\FF_p[v]/(f(v))$, we take each of our 15 elliptic curves $y^2=x^3+A(u)x+B(u)$, with $A,B\in \ZZ[u]$,
and compute $A,B\in \CF$ by reducing $A$ and $B$ modulo $p$ and $f$; this amounts to evaluating them at the image of $v$ in $\CF=\FF_p[v]/(f(v))$.
We now run Schoof's algorithm on $y^2=x^3+Ax+B$ with $A,B\in \CF$, which we can view as running Schoof's algorithm on $d$ distinct elliptic curves over $\FF_p$, one for each root of~$f$.  Unless these~$d$ elliptic curves all have the same trace of Frobenius, the algorithm will eventually test an identity of the form \eqref{eq:schoof} that is satisfied for one but not all of them.
When this happens we will discover a nontrivial zero divisor and a nontrivial factorization of~$f$.  Each root of~$f$ is a generic value of $u\in \FF_p$, and under Conjecture \ref{conj:distincttuples}, at least one of our~15 elliptic curves will yield distinct Frobenius traces when specialized at the roots of~$f$.

In summary, conditional on Conjecture~\ref{conj:distincttuples} we obtain a deterministic polynomial-time algorithm for factoring polynomials over finite fields, as recorded in the theorem below.

\begin{theorem}
Let $f$ be a polynomial over a finite field.  Under Conjecture~\textup{\ref{conj:distincttuples}} there is a deterministic polynomial-time algorithm to compute the factorization of~$f$ into irreducible polynomials.
If $f\in \FF_p[x]$ is the product of $d$ distinct linear polynomials, then a proper divisor of $f$ can be found in $O(d^{1+o(1)}(\log p)^{5+o(1)})$ time.
\end{theorem}

\begin{remark}
Augmenting Kayal's 2005 idea by considering multiple elliptic curves is quite natural, and computationally much more convenient than working with high dimensional abelian varieties as Poonen proposed.
Actually, there is no reason why we couldn't just pick 15 (actually, fewer) random families of elliptic curves and use them in a factoring algorithm.
Heuristically we expect these randomly chosen families to have distinct tuples of Frobenius traces modulo $p$ for all primes $p$, which suffices for the application to factoring.  The downside is that proving this heuristic appears intractable.  

Working with a natural family of abelian varieties is computationally harder, but allows the possibility that geometric methods could prove that the zeta functions over $\FF_p$-specializations are all distinct.
Our proposal bridges the gap, providing an explicit choice of the families of elliptic curves (with inherent computational advantages)
while maintaining a geometric connection that suggests approaches toward proving Conjecture~\ref{conj:distincttuples}.
\end{remark}
\bigskip

\begin{small}
\noindent
\textbf{Data Availability}

\noindent
All data and code associated to this article is available in the GitHub repository \cite{DoublyIsogenousRepo}.
\bigskip

\noindent
\textbf{Conflict of interest}

\noindent
The authors have no conflicts of interest to declare.
\end{small}


\bibliographystyle{hplaindoi}
\bibliography{doubly}

\end{document}